\documentclass[thmsa,final,french,11pt]{article}

\usepackage{eurosym}
\usepackage{amsfonts}
\usepackage{amsmath}
\usepackage{amssymb}
\usepackage{mathrsfs} 
\usepackage{enumerate}
\usepackage{graphicx}

\usepackage{ntheorem}

\newtheorem{Th}{Theorem}
\newtheorem*{Hyp}{Condition (C)}
\newtheorem{Prop}[Th]{Proposition}
\newtheorem{Le}[Th]{Lemma}

\newcommand{\WW}{\mathbf{W}}

\newcommand{\RR}{\mathbf{R}}

\newcommand{\cell}{\mathcal{C}}

\newcommand{\conv}[2][n]{\underset{#1\rightarrow #2}{\longrightarrow}}
\newcommand{\eq}[2][n]{\underset{#1\rightarrow #2}{\sim}}
\newcommand{\EEE}[1]{\operatorname{\mathbb{E}}\left[\,#1\,\right]}

\newcommand{\PPP}[1]{\operatorname{\mathbb{P}}\left(\,#1\,\right)}

\newcommand{\ind}[1]{\mathbb{I}_{#1}\,}

\def\WW_#1{\boldsymbol{W}\!_#1}
\numberwithin{equation}{section}

\newenvironment{prooft}{\noindent {\bf Proof}} {\hfill $\square$  \noindent}

\usepackage{xcolor}

\usepackage[french]{varioref}
\begin{document}

\title{Cluster size distributions of extreme values for the Poisson-Voronoi tessellation}
\author{Nicolas CHENAVIER\footnote{Universit\'e Littoral C\^ote d'Opale, EA 2797, LMPA, 50 rue Ferdinand Buisson, F-62228 Calais, France. E-mail: nicolas.chenavier@univ-littoral.fr}, Christian Y. ROBERT\footnote{Universit\'e de Lyon, Universit\'e Lyon 1, Institut de Science financi\`ere et d'Assurances, 50 avenue Tony Garnier, F-69007 Lyon, France. E-mail: christian.robert@univ-lyon1.fr}}
\maketitle

\begin{abstract}
We consider the Voronoi tessellation based on a homogeneous Poisson point 
process in $\mathbf{R}^{d}$. For a geometric characteristic of the cells
(e.g. the inradius, the circumradius, the volume), we investigate the point
process of the nuclei of the cells with large values. Conditions are
obtained for the convergence in distribution of this point process of
exceedances to a homogeneous compound Poisson point process. We provide a  
characterization of the asymptotic cluster size distribution which is based
on the Palm version of the point process of exceedances. This
characterization allows us to compute efficiently the values of the extremal index and the cluster size probabilities by simulation for various geometric characteristics. The extension to the Poisson-Delaunay
tessellation is also discussed.
\end{abstract}





\section{Introduction}

\label{sec:intro}

\paragraph{Stationary tessellations and the Poisson-Voronoi
tessellation}

A tessellation in $\mathbf{R}^{d}$, $d\geq 1$, endowed with its
Euclidean norm $|\cdot |$, is a countable collection of non-empty convex
compact subsets, called \textit{cells}, with disjoint interiors which
subdivides the space and such that the number of cells intersecting any
bounded subset of $\mathbf{R}^{d}$ is finite. The set $\mathbf{T}$ of
tessellations is endowed with the $\sigma $-field generated by the sets $\{m\in \mathbf{T},\mathbf{\cup }_{C\in m}\partial C\cap K=\varnothing \}$, 
where $\partial K$ is the boundary of $K$ for any compact set $K$ in $\mathbf{R}^{d}$. By a random tessellation $m$, we mean a random variable
with values in $\mathbf{T}$. For a complete account on random
tessellations and their applications, we refer to the books  \cite{SW, SKM}.

A tessellation $m$ is said to be stationary if its distribution is invariant
under translations of the cells. Given a fixed realization of a stationary
tessellation $m$, we associate with each cell $C\in m$, in a deterministic way, a point $z(C)$ which
is called the \textit{nucleus} of the cell, such that $z(C+x)=z(C)+x$ for
all $x\in \mathbf{R}^{d}$. To describe the mean behavior of the
tessellation, the notions of intensity and typical cell are introduced as
follows. Let $A\subset \mathbf{R}^{d}$ be a Borel subset such that $\lambda_d(A)=1$, where  $\lambda _{d}$ is the $d$-dimensional Lebesgue measure. The \textit{intensity} of a stationary tessellation $m$ is defined as 
\begin{equation*}
\gamma_m :=\mathbb{E}\left[ \# \{C\in m,z(C)\in
A\}\right],  
\end{equation*}
where $\# \mathcal{S}$
denotes the cardinality of any finite set $\mathcal{S}$. Thanks to the stationarity of $m$, the intensity does not depend on the
choice of $A$. Without loss of generality, we assume that $\gamma_m =1$.

The \textit{typical cell} $\mathcal{C}$ of a stationary tessellation $m$ is
a random polytope with distribution given by 
\begin{equation}
\mathbb{E}[f(\mathcal{C})]=\mathbb{E}\left[
\,\sum_{C\in m,z(C)\in A}f(C-z(C))\,\right] ,  \label{def:typicalcell}
\end{equation}
where  $f:\mathcal{K}_{d}\rightarrow \mathbf{R}$
is any bounded measurable function on the set of convex bodies $\mathcal{K}_{d}$ (endowed with the Hausdorff topology). 

Let $\chi $ be a locally finite subset of $\mathbf{R}^{d}$. The Voronoi
cell with nucleus $x\in \chi $ is the set of all sites $y\in \mathbf{R}^{d}$
whose distance from $x$ is smaller or equal than the distances to all other
points of $\chi $, i.e. 
\begin{equation*}
C_{\chi }(x):=\{y\in \mathbf{R}^{d}:|y-x|\leq |y-x^{\prime }|,x^{\prime }\in
\chi \}.
\end{equation*}
When $\chi =\eta $ is a homogeneous Poisson point process, the family $
m:=\{C_{\eta }(x):x\in \eta \}$ is the so-called \textit{Poisson-Voronoi
tessellation}. The intensity of such a tessellation equals the intensity of $\eta$. A consequence of the theorem of Slivnyak (see e.g.
Theorem 3.3.5 in \cite{SW}) shows that 
\begin{equation}
\label{eq:typicalcell}
\mathcal{C}\overset{\mathcal{D}}{=}C_{\eta \cup \{0\}}(0),
\end{equation}
where $\overset{\mathcal{D}}{=}$ denotes the equality in distribution. The study of this typical cell in the literature includes mean values
calculations \cite{M}, second order properties \cite{HM} and distributional
estimates \cite{BL,Mu}. Voronoi tessellations are extensively used in many domains such as cellular biology \cite{Po}, astrophysics \cite{Za}, telecommunications \cite{BB3} and finance \cite{PW}. For a complete account on Poisson-Voronoi tessellations and their applications, we refer to the book by Okabe et \textit{al.} (see Chapter 5 in \cite{OBSC}). 

\paragraph{Point process of exceedances for a stationary sequence of
real random variables}

Let $\left( X_{n}\right) _{n\in \mathbb{Z}}$ be a strictly stationary
sequence of real random variables. Assume
that for each $\tau >0$ there exists a sequence of levels $\left(
u_{n}\left( \tau \right) \right) $ such that $\lim_{n\rightarrow \infty }n\mathbb{P}\left( X_{1}>u_{n}\left( \tau \right) \right) =\tau $. The point
process of time normalized exceedances is defined by $\phi _{B}(\tau):=n^{-1}\cdot
\left\{ i\in B:X_{i}>u_{n}\left( \tau \right) \right\} $ for any Borel set $B\subset W_{n}:=[-n/2,n/2]$. If $\left( X_{n}\right) $ satisfies a long
range dependence condition (known as condition $\Delta (u_{n}(\tau ))$) and if 
the point process $\phi _{W_{n}}(\tau)$ weakly converges to a point process in $[-1/2,1/2]$, then the limiting point process is necessarily a homogeneous
compound Poisson process with intensity $\nu\geq 0 $ and limiting
\textit{cluster size distribution} $\pi $ (see Corollary 3.3 in \cite{HHL}). According to Leadbetter \cite{L2}, the
constant $\theta =\nu /\tau $ is referred to as the \textit{extremal index}. It may
be shown that $0\leq \theta \leq 1$ and that the compound Poisson limit becomes
Poisson when $\theta =1$.\label{page:extremalindex}

If $\lim_{n\rightarrow \infty }\mathbb{P}(\# \phi _{W_{n}}(\tau)=0)=e^{-\theta
\tau }$, then a necessary and sufficient condition for the convergence of $\phi
_{W_{n}}$ is the convergence of the conditional distribution of $\# \phi
_{B_{n}}$, with $B_{n}=\left[ 0,q_{n}\right]$, given that there is at least
one exceedance of $u_{n}(\tau )$ among $X_{1},\ldots ,X_{q_{n}}$, to the
distribution $\pi =\left( \pi _{k}\right) _{k\geq 1}$, i.e. 
\begin{equation}
\lim_{n\rightarrow \infty }\mathbb{P}\left( \left. \# \phi
_{B_{n}}(\tau)=k\right\vert \# \phi _{B_{n}}(\tau)>0\right) :=\pi _{k},\text{\qquad }
k\geq 1,  \label{condcv}
\end{equation}
where $\left( q_{n}\right) $ is a $\Delta (u_{n}(\tau ))$-separating
sequence, with $\lim_{n\rightarrow \infty }q_{n}/n=0$ (see Theorem 4.2 in \cite{HHL}). This condition is known as the \textit{blocks}
characterization of the cluster size distribution $\pi$. Under additional mild
conditions (see e.g. \cite{Sm}) the extremal index is equal to the
reciprocal of the mean of $\pi $.

An equivalent condition to $\eqref{condcv}$ is proposed in Theorem 4.1\ in 
\cite{Roo} (see also Theorem 2.5 in \cite{Pe}) and is given by
\begin{equation}
\lim_{n\rightarrow \infty }\mathbb{P}\left( \left. \# \phi
_{B_{n}}(\tau)=k\right\vert X_{0}>u_{n}\left( \tau \right) \right) :=p'_{k}=\theta
\sum_{m=k}^{\infty }\pi _{m},\text{\qquad }k\geq 1.  \label{condcv2}
\end{equation}
In particular, we have $\theta =p'_{1}$. This second condition is useful to compute
the values of the extremal index and the cluster size probabilities when the
conditional distributions of the exceedances may be derived from the
dynamics of $\left( X_{n}\right) _{n\in \mathbb{Z}}$, e.g. for the
regularly varying multivariate time series \cite{BS} or the Markov
sequences \cite{Pe}. This condition may be called the \textit{runs}
characterization of the cluster size distribution since the runs estimator
of the extremal index is based on the following result $\theta
=\lim_{n\rightarrow \infty }\mathbb{P}\left( \left. \cap
_{i=1}^{q_{n}}\{X_{i}\leq u_{n}\left( \tau \right) \}\right\vert
X_{0}>u_{n}\left( \tau \right) \right)$. The runs characterization is natural for a random object
as a time series where the direction of time is used to design the dynamics
of the series. Estimators of the extremal index and the cluster size distribution, based on the blocks and  runs characterizations, are extensively investigated, see e.g. \cite{Rob,S2}.

However, we claim that it could also be useful to consider a new condition
where the conditional event $\{X_{0}>u_{n}\left( \tau \right) \}$ is not
used as the starting point of the considered cluster, but as a part of this
cluster. We therefore introduce a new discrete probability distribution $p=\left( p_{k}\right) _{k\geq 1}$ and the following condition 
\begin{equation}
\lim_{n\rightarrow \infty }\mathbb{P}\left( \left. \# \phi
_{C_{n}}(\tau)=k\right\vert X_{0}>u_{n}\left( \tau \right) \right) :=p_{k},\text{\qquad }k\geq 1,  \label{condcv3}
\end{equation}
where $C_{n}=\left[ -q_{n}/2,q_{n}/2\right] $. If $p$ exists, an adaptation of our main result (see Theorem \ref{Th:theta}) shows that $p_{k}=\theta k\pi _{k}$ for $k\geq 1$, and therefore $\theta=\sum_{k=1}^{\infty }k^{-1}p_{k}$. Such a condition will be proposed for
random tessellations for which there is no natural direction in the space $\mathbf{R}^{d}$. However, we think that our new condition could be fruitful for time series.

\paragraph{Point process of exceedances for a stationary tessellation}

Let $m$ be a stationary tessellation in $\mathbf{R}^{d}$. We consider a
geometric characteristic $g:\mathcal{K}_{d}\rightarrow \mathbb{R}$, which is a measurable translation-invariant function, i.e. $g(C+x)=g(C)$ for all $C\in \mathcal{K}_{d}$ and $x\in \mathbf{R}^{d}$, and such that, for any $\tau >0$, there exists a threshold $v_{\rho
}(\tau )$ satisfying 
\begin{equation}
\lim_{\rho \rightarrow \infty }\rho \mathbb{P}\left( g(\mathcal{C})>v_{\rho
}(\tau )\right) =\tau,  \label{Threshcond}
\end{equation}
where $\mathcal{C}$ is the typical cell. We observe only a part of the stationary tessellation $m$ in the window $\mathbf{W}_{\rho }:=\rho ^{1/d}\cdot \left[ -1/2,1/2\right] ^{d}$, $\rho >0$,
and we are interested in the point process of exceedances $\Phi _{\mathbf{W}_{\rho }}(\tau)$ where, for any Borel set $B\subset \mathbf{R}^{d}$, we let  
\begin{equation*}
\Phi _{B}(\tau):=\rho ^{-1/d}\cdot \left\{ z(C):z(C)\in B, g(C)>v_{\rho }(\tau
), C\in m\right\} .
\end{equation*}
In this paper, we investigate the weak convergence of the point process $\Phi _{\mathbf{W}_{\rho }}(\tau)$ in $\left[ -1/2,1/2\right] ^{d}$ as $\rho $ tends to infinity. In 
\cite{Chen}, a first result was obtained for geometric characteristics for
which a short range dependence condition holds (equivalent to the  so-called condition  $D^{\prime }$ for stationary sequences of real random variables): it is shown
that the point process $\Phi _{\mathbf{W}_{\rho }}(\tau)$ weakly converges to a 
homogeneous Poisson point process with intensity $\tau $. In this paper, we
are interested in finding weaker conditions for other geometric
characteristics such that the point process $\Phi _{\mathbf{W}_{\rho }}(\tau)$
weakly converges to a homogeneous compound Poisson point process.

Let $\mathfrak{B}_{\rho}$ be a sub-cube of $\mathbf{W}_{\rho }$ such
that $\lim_{\rho \rightarrow \infty }\lambda _{d}(\mathfrak{B}_{\rho})/\rho =0$. Condition $\eqref{condcv}$ for the tessellation $m$ will be written in the
following way: 
\begin{equation}
\lim_{\rho \rightarrow \infty }\mathbb{P}\left( \left. \# \Phi _{\mathfrak{B}_{\rho}}(\tau)=k\right\vert \# \Phi _{\mathfrak{B}_{\rho}}(\tau)>0\right) =\pi
_{k},\text{\qquad }k\geq 1,   \label{Blockcond}
\end{equation}
for a discrete probability distribution $\pi =\left( \pi _{k}\right) _{k\geq
1}$, which we also call the cluster size distribution. Additional assumptions on $\mathfrak{B}_{\rho}$ will be necessary and will
depend on the mixing properties of the tessellation. Condition $\eqref{condcv2}$ cannot be transposed for stationary tessellations as
explained previously. Condition $\eqref{condcv3}$ has to be modified since
the cell which contains the origin (the Crofton cell) is not distributed as
the  typical cell. To overcome this difficulty, we consider a Palm version $\Phi_{\RR^d} ^{0}(\tau)$ of $\Phi_{\RR^d}(\tau) $, i.e. a point process whose distribution is given by the
Palm distribution of $\Phi_{\RR^d}(\tau)$ (see Sections 3.3 and 3.4 in \cite{SW} for a complete account on Palm theory). For any $B\subset \RR^d$, we also let $\Phi _{B}^{0}(\tau)=\Phi_{\RR^d}
^{0}(\tau)\cap B$ . An analogous version of Condition $\eqref{condcv3}$ in the context of random tessellations can be stated as follows: 
\begin{equation}
\lim_{\rho \rightarrow \infty }\mathbb{P}\left( \# \Phi _{\mathfrak{B}
_{\rho}}^{0}(\tau)=k\right) :=p_{k},\text{\qquad }k\geq 1,  \label{Palmcond}
\end{equation}
for a discrete probability distribution $p=\left( p_{k}\right) _{k\geq 1}$.
\ In general the distributions $\pi $ and $p$ cannot be made explicit. It is necessary to use simulations to compute approximate values  of the probabilities $\pi _{k}$ and $p_{k}$. 

The blocks method $\eqref{Blockcond}$ competes with the Palm approach $\eqref{Palmcond}$. The idea of the Palm approach is to consider clusters
close to the origin given that the cell whose nucleus is the origin has an
exceedance. Our approach provides better approximations of the extremal index and the cluster size distribution and requires less simulations. Indeed, we simulate only blocks that contain at least one exceedance (the one of the Crofton cell that contains
the origin), while with the blocks approach, it is necessary to simulate a
very large number of blocks (including those without any extreme value). More precisely, in our numerical illustrations in $\RR^2$, we simulate tessellations only observed in the square $[-173, 173]^2$ to approximate $\theta$ and $p=(p_k)_{k\geq 1}$  thanks to our Palm approach. A blocks approach would have required to simulate tessellations in the square $[-5.18\cdot 10^{21},5.18\cdot 10^{21} ]$, which is practically impossible.  

In this paper, we only focus on the Poisson-Voronoi tessellation and we 
 explain how our results can be extended to the 
Poisson-Delaunay tessellation. In Section \ref{sec:preliminaries}, we give several preliminaries by introducing notation and conditions on our geometric characteristic. In Section \ref{sec:weakconvergence}, we investigate the convergence in distribution of the point
process of exceedances to a homogeneous compound Poisson point process. This convergence is stated in our main result (Theorem \ref{Th:theta}).  In Section \ref{sec:examples}, we
give three examples and numerical illustrations. The extension to 
the Poisson-Delaunay tessellation is discussed in Section \ref{sec:Delaunay}. 

\section{Preliminaries}
\label{sec:preliminaries}
In this section, we introduce several notation and conditions which will be used throughout the paper.

\paragraph{Notation}

\begin{itemize}
\item Let $x\in \mathbf{R}^{d}$ and let $A,B\subset \mathbf{R}^{d}$ be two subsets. We write $x+A:=\{x+a:a\in A\}$, $A\oplus B:=\{a+b:a\in A,b\in
B\}$ and $A\ominus B:=\{x\in \mathbf{R}^{d}:x+B\subset A\}$. Moreover, we denote the complement of $A$ by $A^c:=\RR^d\setminus A$.

\item For any $A,B\subset\mathbf{R}^d$, we denote the distance between $A$ and $B$ by 
 $\delta(A,B):=\inf_{(a,b)\in A\times B}|a-b|$.

\item For any $k$-tuple of points $x_{1},\ldots ,x_{k}\in \mathbf{R}^{d}$,
we write $x_{1:k}:=(x_{1},\ldots ,x_{k})$. With a slight abuse of notation, we also write $\{x_{1:k}\}:=\{x_{1},\ldots ,x_{k}\}$. 

\item For any Borel subset $B\subset \mathbf{R}^{d}$, we write $\mathbf{B}_\rho:=\rho ^{1/d}\cdot B$. 

\item We denote by $\mathcal{F}_{lf}$ the set of locally finite subsets in $\RR^d$. This set is endowed with the $\sigma$-field induced by the so-called Fell topology on $\mathcal{F}_{lf}$ (see e.g. p. 563 in \cite{SW}). 

\item Let $\chi\in \mathcal{F}_{lf}$. 

\begin{itemize}
\item For any $x_{1:k}\in \chi ^{k}$ and for any $v\geq 0$, we write $g^{\chi }(x_{1:k})>v$ to specify that $g(C_{\chi }(x_{j}))>v$ for any $1\leq
j\leq k$. In particular, we let $g^\chi(x):=g(C_\chi(x))$. 

\item For any $B\subset \mathbf{R}^{d}$, we write $M_{B}^{\chi }:=\max_{x\in
\chi \cap B}g^{\chi }(x)$. When $\chi \cap B=\varnothing $, we take $M_{B}^{\chi }:=-\infty $.

\item For any $\rho >0$ and $\tau >0$ , we denote by  $\Phi^{\chi }(\tau)$ the point process of exceedances, i.e. 
\begin{equation*}
\Phi^{\chi }(\tau):=\rho ^{-1/d}\cdot \left\{ x\in \chi: g^{\chi
}(x)>v_{\rho }(\tau )\right\}.
\end{equation*}
Besides, for any $B\subset \mathbf{R}^{d}$, we write $\Phi_B^{\chi }(\tau):=\Phi^\chi(\tau) \cap (\rho^{-1/d}B)$.
\end{itemize}

\item We denote by $\eta$ a homogeneous Poisson point process in $\RR^d$. Excepted in Section \ref{sec:Delaunay}, we assume that the intensity of $\eta$ is $\gamma_\eta=1$.

\item For each $\tau>0$, we denote by $\Phi^{\eta, 0}(\tau)$ the Palm version of $\Phi^{\eta}(\tau)$. In particular, for any $B\subset\RR^d$ we let $\Phi _{B}^{\eta, 0 }(\tau):=\Phi^{\eta, 0 }(\tau)\cap (\rho^{-1/d} B)$. We also associate two probabilities defined as follows:
\begin{equation*}
\pi _{k,B}(\tau):=\mathbb{P}\left( \#\Phi_{B}^{\eta}(\tau)=k|\#\Phi_{B}^{\eta}(\tau)>0\,\right) \quad 
\text{ and }\quad p_{k,B}(\tau):=\mathbb{P}(\#\Phi_{B}^{\eta, 0}(\tau)=k\,).
\end{equation*}
The quantity $\pi _{k,B}(\tau)$ is the probability that there are $k$ exceedances in $B$ conditional on the fact that there is at less an exceedance in $B$, whereas the quantity $p_{k,B}(\tau)$ is the probability that there are $k$ exceedances in $B$ conditional on the fact that the origin is a nucleus and that the cell with nucleus the origin is an exceedance. Notice that these probabilities  also depend on $\rho$. 

\item For any pair of functions $h_{1},h_{2}\colon \mathbf{R}\rightarrow 
\mathbf{R}$, we write $h_{1}(\rho )\underset{\rho \rightarrow \infty }{\sim }h_{2}(\rho )$ and $h_{1}(\rho )=O(h_{2}(\rho ))$ to respectively mean that $h_{1}(\rho )/h_{2}(\rho )\rightarrow 1$ as $\rho \rightarrow \infty $ and $h_{1}(\rho )/h_{2}(\rho )$ is bounded for $\rho $ large enough.

\item We denote by $q: \rho\mapsto q_{\rho }$ a generic function such that, for any $\alpha ,\beta
>0$, we have simultaneously 
\begin{equation}
q_{\rho }\cdot (\log \rho )^{\alpha }\cdot \rho ^{-1}\underset{\rho
\rightarrow \infty }{\longrightarrow }0\quad \text{ and }\quad q_{\rho
}^{-1}\cdot (\log \rho )^{\beta }\underset{\rho \rightarrow \infty }{\longrightarrow }0.  \label{def:qrho}
\end{equation}

\item Let $\varepsilon>0$ be fixed. For any $\rho>0$, we denote by $n_{\rho}$ and $m_{\rho}$ the integers \begin{equation}\label{def:nm} n_{\rho}:=\left\lfloor (\log \rho )^{-(1+\varepsilon)/d}\cdot \rho ^{1/d}\right\rfloor \quad \text{ and } \quad m_{\rho}:=\left\lfloor  q_{\rho }^{-1/d}\cdot (\log \rho )^{-(1+\varepsilon )/d}\cdot \rho ^{1/d}\right\rfloor. 
\end{equation} We also define two squares centered at $0$ as follows:
\begin{equation}
\mathfrak{c}_{\rho}:= \frac{\rho^{1/d}}{n_{\rho}}\cdot \lbrack -D,D]^{d} \quad\text{ and }\quad   \mathfrak{Q}_{\rho}:= \frac{\rho^{1/d}}{m_{\rho}}\cdot \lbrack -1/2,1/2]^{d},  \label{def:B0}
 \end{equation}
 where $D:=2\cdot (\lfloor \sqrt{d}\rfloor +1)$. In particular, we have \[\lambda_d(\mathfrak{c}_{\rho})\eq[\rho]{\infty}(\log\rho)^{1+\varepsilon}(2D)^d \quad\text{and}\quad  \lambda_d(\mathfrak{Q}_{\rho}) \eq[\rho]{\infty}q_\rho(\log\rho)^{1+\varepsilon}.\]
\end{itemize}

Throughout the paper, we use $c$ to signify a universal positive constant
not depending on $\rho $ but which may depend on other quantities. When
required, we assume that $\rho $ is sufficiently large.

\paragraph{Conditional independence}

Let $\eta $ be a homogeneous Poisson point process. We begin with a first lemma that characterizes the dependence structure of the Poisson-Voronoi tessellation induced by $\eta$. 

We partition $\mathbf{W}_\rho$ into a set $\mathcal{V}_{\rho}$ of $n_{\rho}^{d}$ sub-cubes of equal size, where  $n_{\rho}^{d}$ is defined in \eqref{def:nm}. These sub-cubes are indexed
by the set of $\mathfrak{i}:=(i_{1},\ldots ,i_{d})\in \left[ 1,n_{\rho}\right] ^{d}$. With a slight abuse of notation, we identify a cube with its index. Notice that for each $\mathfrak{i}\in \mathcal{V}_{\rho}$, we have $\lambda_d(\mathfrak{i}) = (2D)^{-d}\cdot \lambda_d(\mathfrak{c}_{\rho})$. Besides, the distance between sub-cubes $\mathfrak{i}$ and $\mathfrak{j}$ is denoted by $d(\mathfrak{i},\mathfrak{j}):=\max_{1\leq r\leq d}|i_{r}-j_{r}|$. Moreover, if $\mathcal{A}$ and $\mathcal{B}$ are two sets of sub-cubes, we let $d(\mathcal{A},\mathcal{B})=\min_{\mathfrak{i}\in \mathcal{A},\mathfrak{j}\in \mathcal{B}}d(\mathfrak{i},\mathfrak{j})$ and 
\begin{equation}
\label{eq:notalgebra}
\Sigma _{\mathcal{A}}^{\eta \cup \{x_{1:k}\}}:=\sigma \{g^{\eta \cup
\{x_{1:k}\}}(x):x\in (\eta \cup \{x_{1:k}\})\cap \mathfrak{i},\mathfrak{i}%
\in \mathcal{A}\}.
\end{equation}
Finally, to ensure several independence properties, we introduce the following event: 
\begin{equation*}
\mathscr{A}_{\rho }:=\cap _{\mathfrak{i}\in \mathcal{V}_{\rho}}\{\eta \cap \mathfrak{i}\neq
\varnothing \}.
\end{equation*}
The event $\mathscr{A}_{\rho }$ is extensively used in stochastic geometry to derive central
limit theorems or to deal with extremes (see e.g. \cite{AB,Chen}). It will
 play a crucial role in the rest of the paper. The following lemma is the heart of our development and captures the idea of ``local dependence''.  
\begin{Le}
\label{Le:Arho} Let  $x_{1},\ldots ,x_{k}\in \mathbf{R}^{d}$, with $k\geq 0$. Then

\begin{enumerate}[(i)]
\item \label{Le:Arho1} conditional on $\mathscr{A}_{\rho }$, the $\sigma $-fields $\Sigma _{\mathcal{A}}^{\eta \cup \{x_{1:k}\}}$ and $\Sigma _{\mathcal{B}}^{\eta \cup
\{x_{1:k}\}} $ are independent when $d(\mathcal{A},\mathcal{B})>D$;

\item \label{Le:Arho2} for any $\alpha >0$, we have $\rho ^{\alpha }\cdot \mathbb{P}\left(
\,\mathscr{A}_{\rho}^{c}\,\right) \underset{\rho \rightarrow \infty }{\longrightarrow 
}0$.
\end{enumerate}
\end{Le}

\begin{prooft}
The first assertion is a simple adaptation of Lemma 5 in \cite{Chen}. The
second one comes from \eqref{def:nm} and the fact that 
\begin{equation*}
\mathbb{P}\left( \,\mathscr{A}_{\rho }^{c}\,\right) =\mathbb{P}\left( \,\bigcup_{\mathfrak{i}\in \mathcal{V}_{\rho}}\{\eta \cap \mathfrak{i}=\varnothing \}\,\right)
\leq n_{\rho}^{d}e^{-\rho /n_{\rho}^{d}}.
\end{equation*}
\end{prooft}

\paragraph{Condition on the geometric characteristic}
To state our main theorem, we assume some condition on the geometric characteristic $g$, referred to as Condition (C). 
\begin{Hyp}
For any $\tau>0$, there exists a constant $c$ such that, for
any $(k-1)$-tuple of points $y_{2:k}\in \mathbf{R}^{d(k-1)}$, for any $z\in \{0,y_{2:k}\}$, we have 
\begin{equation*}
\mathbb{P}\left( \,g^{\eta \cup \{0,y_{2:k}\}}(z)>v_{\rho }(\tau)\,\right) \leq
c\cdot \rho ^{-1},
\end{equation*}
where $v_{\rho }(\tau )$ satisfies Equation $\eqref{Threshcond}$,
with the convention $\{y_{2:k}\}=\varnothing $ when $k=1$.
\end{Hyp}

In particular, Condition (C) is satisfied when $g^{\eta \cup
\{0,y_{2:k}\}}(x)\leq g^{\eta \cup \{0\}}(x)$ for any $x\in \eta \cup \{0\}$
and for any $y_{2:k}\in \mathbf{R}^{d(k-1)}$, i.e. when the geometric
characteristic of a cell always decreases if new points are added to the
point process $\eta$.

\section{Weak convergence of the point process of exceedances for a Poisson-Voronoi tessellation}
\label{sec:weakconvergence}

\subsection{An explicit representation for $p_{k,B}$}
According to the theorem of Slivnyak, the Palm distribution of $\eta $ is given by
the distribution of $\eta \cup \{0\}$. As a consequence, the following lemma
shows that for any $B\subset \mathbf{R}^{d}$, the distribution of $\#\Phi_{B}^{\eta, 0 }(\tau)$ is the same as the one of $\#\Phi_{B}^{\eta \cup \{0\}}(\tau)$ given that $g^{\eta \cup \{0\}}(0)>v_{\rho }(\tau )$.

\begin{Le}
\label{Le:PalmVoronoi}
For any $B\subset \mathbf{R}^{d}$, $\rho >0$ and $k\geq 1$, we have 
\begin{equation*}
p_{k,B}(\tau)=\mathbb{P}(\#\Phi_B^{\eta \cup \{0\}}(\tau)=k|g^{\eta \cup \{0\}}(0)>v_{\rho
}(\tau )).
\end{equation*}
\end{Le}

\begin{prooft}
Since $p_{k,B}(\tau):=\PPP{\#\Phi_B^{\eta, 0} (\tau)= k} = \PPP{\# \Phi^{\eta, 0}(\tau)\cap (\rho^{-1/d}B) = k}$, we obtain for any Borel subset $A\subset \RR^d$, with  $\lambda_d(A)=1$, that  
\begin{equation}\label{eq:Palm1}p_{k,B}(\tau) = \frac{1}{\gamma_{\Phi^\eta(\tau)}}\EEE{ \,\sum_{z\in \Phi ^{\eta }(\tau)\cap A}\mathbb{I}_{\# (\Phi^{\eta }(\tau)-z)\cap
(\rho ^{-1/d}B)=k}\,},   \end{equation}
where $\gamma_{\Phi^\eta(\tau)}:= \EEE{\# (\Phi^\eta(\tau) \cap A)}$ is the intensity of $\Phi^\eta(\tau)$. According to \eqref{def:typicalcell}, this intensity equals
\begin{equation}\label{eq:Palm2}\gamma_{\Phi^\eta(\tau)}  = \EEE{\sum_{x\in \eta\cap \mathbf{A}_\rho} \ind{g^\eta(x)>v_\rho(\tau)}  } = \rho\PPP{g(\cell)>v_\rho(\tau)}. 
\end{equation} Moreover, it results from the Slivnyak-Mecke formula (e.g. Corollary 3.2.3 in \cite{SW}) that for any $B\subset\RR^d$, 
\begin{multline*}
\EEE{ \,\sum_{z\in \Phi ^{\eta }(\tau)\cap A}\mathbb{I}_{\# (\Phi^{\eta }(\tau)-z)\cap
(\rho ^{-1/d}B)=k}\,}\\
\begin{split}
& = \EEE{\sum_{x\in \eta\cap \mathbf{A}_\rho}\ind{\#\{y\in (\eta-x)\cap B: g^\eta(y)>v_\rho(\tau)\} = k}\ind{g^\eta(x)>v_\rho(\tau)}    }\\
& = \int_{\mathbf{A}_\rho}\PPP{\#\{y\in (\eta\cup\{x\}-x)\cap B: g^{\eta\cup\{x\}}(y)>v_\rho(\tau)\} = k, g^{\eta\cup\{x\}}(x)>v_\rho(\tau)}\mathrm{d}x.
\end{split}
\end{multline*}
Thanks to the stationarity of $\eta$ and because $g$ is translation-invariant, the above integrand does not depend on $x$. By integrating over $x\in \mathbf{A}_\rho$ and using the fact that $\lambda_d(\mathbf{A}_\rho) = \rho$, it follows that
\begin{multline*}
\EEE{ \,\sum_{z\in \Phi ^{\eta }(\tau)\cap A}\mathbb{I}_{\# (\Phi^{\eta }(\tau)-z)\cap
(\rho ^{-1/d}B)=k}\,}\\
\begin{split}
& = \rho \cdot \PPP{\#\{y\in (\eta\cup\{0\})\cap B: g^{\eta\cup\{0\}}(y)>v_\rho(\tau)\} = k, g^{\eta\cup\{0\}}(0)>v_\rho(\tau)}\\
& = \rho\cdot \PPP{\#\Phi_B^{\eta\cup\{0\}}(\tau) = k,g^{\eta\cup\{0\}}(0)>v_\rho(\tau) }. 
\end{split}
\end{multline*}
This together with \eqref{eq:typicalcell}, \eqref{eq:Palm1} and \eqref{eq:Palm2} concludes the proof of Lemma \ref{Le:PalmVoronoi}. 
\end{prooft}

\subsection{A technical result}
 The following technical proposition will be the key ingredient to prove our main theorem. 

\begin{Prop}
\label{Prop:identification} Assume that $g$ satisfies Condition (C). Then, for any 
$k\geq 1$, we have 
\begin{equation}
k\mathbb{P}\left( \,\#\Phi_{\mathfrak{Q}_{\rho}}^{\eta }(\tau )=k\,\right) -\lambda
_{d}(\mathfrak{Q}_{\rho})\cdot \mathbb{P}\left( \,\#\Phi_{\mathfrak{Q}_{\rho}}^{\eta
\cup \{0\}}(\tau )=k,g^{\eta \cup \{0\}}(0)>v_{\rho }(\tau )\,\right)
=o\left( \lambda _{d}(\mathfrak{Q}_{\rho})\cdot \rho ^{-1}\right).
\label{eq:identification}
\end{equation}
\end{Prop}

The previous proposition is obvious if we replace $o\left( \lambda
_{d}(\mathfrak{Q}_{\rho})\cdot \rho ^{-1}\right) $ by $O\left( \lambda _{d}(\mathfrak{Q}_{\rho})\cdot \rho ^{-1}\right) $ in \eqref{eq:identification} since $\mathbb{P}\left( \,\#\Phi_{\mathfrak{Q}_{\rho}}^{\eta }(\tau )=k\,\right)=O\left( \lambda _{d}(\mathfrak{Q}_{\rho})\cdot \rho ^{-1}\right) $ and $\mathbb{P}\left( \,\#\Phi_{\mathfrak{Q}_{\rho}}^{\eta
\cup \{0\}}(\tau )=k,g^{\eta \cup \{0\}}(0)>v_{\rho }(\tau )\,\right)=O\left(\rho^{-1}\right)$. Actually, the main difficulty is to prove that the left-hand side is negligible compared
to $\lambda _{d}(\mathfrak{Q}_{\rho})\cdot \rho ^{-1}$, which constitutes the
main ingredient to prove Theorem \ref{Th:theta}.

Since $q_{\rho }$ is any function such that \eqref{def:qrho} holds, we can take $q_{\rho }=(\log \log \rho )^{\log \log \rho }$. Actually, we think that Proposition \ref{Prop:identification} remains true
when $q_{\rho }=\log \rho $, which is slightly more efficient for simulating estimators of $\theta $ and $p _{k}$.


\begin{prooft}
We begin with the case $k\geq 2$. First, we give an integral representation of the left-hand side of \eqref{eq:identification}. Because of the stationarity of $\eta $ and thanks to the Slivnyak-Mecke formula, we have 
\begin{align*}
\mathbb{P}\left( \,\#\Phi_{\mathfrak{Q}_{\rho}}^{\eta }(\tau)=k\,\right) & =\frac{1}{k!}
\mathbb{E}\left[ \,\sum_{x_{1:k}\in (\eta \cap \mathfrak{Q}_{\rho})^{k}}\mathbb{I}_{g^{\eta }(x_{1:k})>v_{\rho}(\tau)}\,\mathbb{I}_{M_{\mathfrak{Q}_{\rho}\setminus
\{x_{1:k}\}}^{\eta }\leq v_{\rho}(\tau)}\,\,\right]  \\
& =\frac{1}{k!}\int_{\mathfrak{Q}_{\rho}}\int_{(\mathfrak{Q}_{\rho}-x_{1})^{k-1}}p_{x_{1}}(y_{2:k})\mathrm{d}y_{2:k}\mathrm{d}x_{1},
\end{align*}
where, for any $x_{1}\in \mathbf{R}^{d}$ and any $y_{2:k}\in \mathbf{R}^{(k-1)d}$, we write
\begin{equation*}
p_{x_{1}}(y_{2:k}):=\mathbb{P}\left( \,g^{\eta \cup
\{0,y_{2:k}\}}(0,y_{2:k})>v_{\rho}(\tau),M_{(\mathfrak{Q}_{\rho}-x_{1})\setminus
\{0,y_{2:k}\}}^{\eta \cup \{0,y_{2:k}\}}\leq v_{\rho}(\tau)\,\right) .
\end{equation*}
In the same spirit as above, we also obtain 
\begin{equation*}
\mathbb{P}\left( \,\#\Phi_{\mathfrak{Q}_{\rho}}^{\eta \cup \{0\}}(\tau)=k,g^{\eta \cup
\{0\}}(0)>v_{\rho}(\tau)\,\right) =\frac{1}{(k-1)!}\int_{\mathfrak{Q}_{\rho}^{k-1}}p_{0}(y_{2:k})\mathrm{d}y_{2:k}.
\end{equation*}
Integrating the right-hand side of the above equation over $\mathfrak{Q}_{\rho}$, it follows that 
\begin{multline*}
k\mathbb{P}\left( \,\#\Phi_{\mathfrak{Q}_{\rho}}^{\eta }(\tau)=k\,\right) -\lambda _{d}(\mathfrak{Q}_{\rho})\cdot \mathbb{P}\left( \,\#\Phi_{\mathfrak{Q}_{\rho}}^{\eta \cup
\{0\}}(\tau)=k,g^{\eta\cup\{0\} }(0)>v_{\rho}(\tau)\,\right)  \\
=\frac{1}{(k-1)!}\int_{\mathfrak{Q}_{\rho}}\left( \int_{(\mathfrak{Q}_{\rho}-x_{1})^{k-1}}p_{x_{1}}(y_{2:k})\mathrm{d}y_{2:k}-\int_{\mathfrak{Q}_{\rho}^{k-1}}p_{0}(y_{2:k})\mathrm{d}y_{2:k}\right) \mathrm{d}x_{1}.
\end{multline*}

The main difficulty to prove that the right-hand side equals $o\left( \lambda _{d}(\mathfrak{Q}_{\rho})\cdot \rho ^{-1}\right)$ comes from the dependence between the $(k+1)$-events considered in the probability $p_{x}(y_{2:k})$, with $x=x_1$ and $x=0$. Actually, the more the distances between the $y_{2},\ldots ,y_{k}$ is large, the more the dependence is weak. To overcome this difficulty,
we introduce the following event:  
\begin{equation*}
E_{m}:=\left\{ y_{2:k}\in \mathbf{R}^{(k-1)d}:S_0(y_{2:k})\text{ has $m$
connected components }\right\} ,
\end{equation*}
where, for any $y_{2:k}\in \mathbf{R}^{d(k-1)}$, the set $S_0(y_{2:k})\subset 
\mathbf{R}^{d}$ is defined as (see Figure \ref{fig:Figset})
\begin{equation*}
S_0(y_{2:k}):=\mathfrak{c}_{\rho}\cup \bigcup_{j=2}^{k}(y_{j}+\mathfrak{c}_{\rho}).
\end{equation*}
\begin{center}
\begin{figure}
\begin{center}
\begin{tabular}{cccc}
    \includegraphics[width=7.5cm,height=7.5cm]{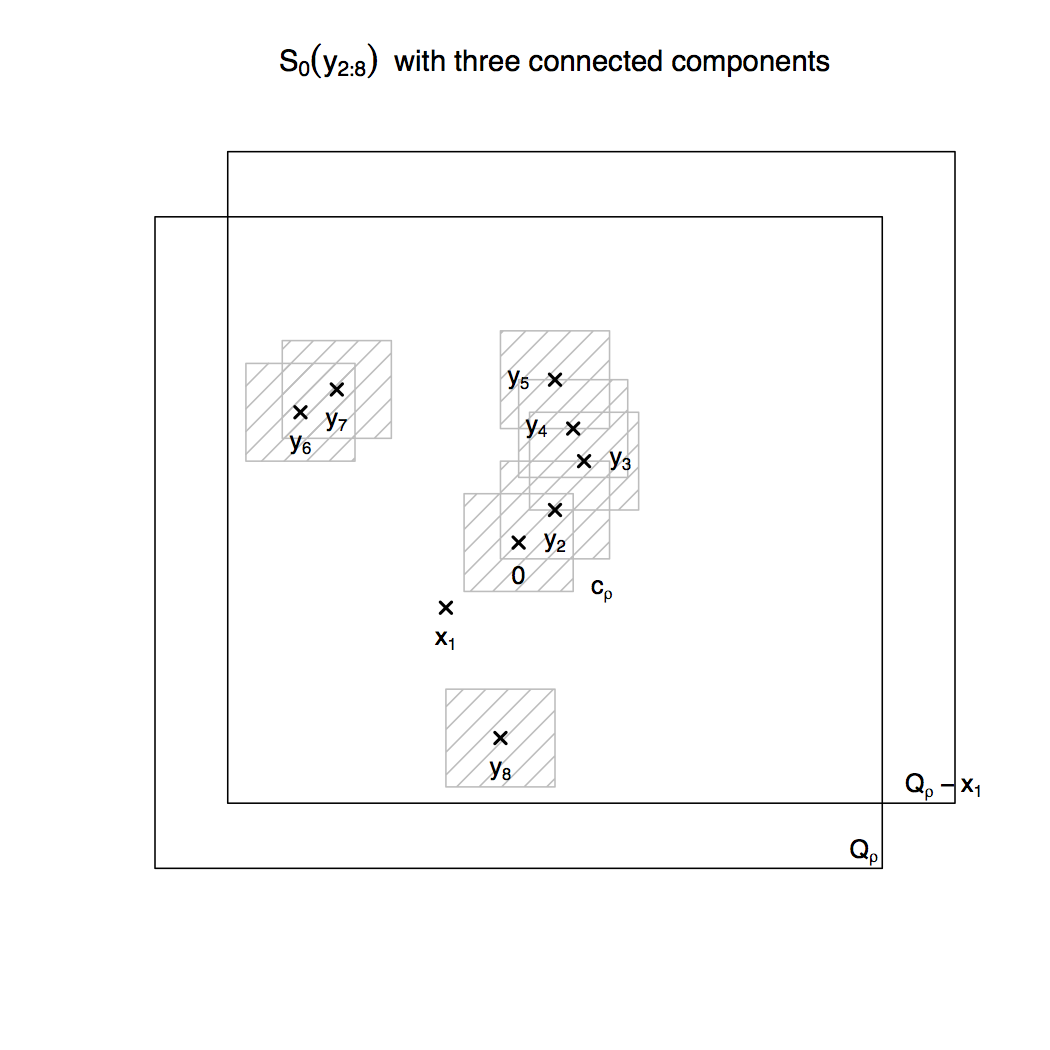} 
\end{tabular}
\end{center}
\caption{\label{fig:Figset} A configuration of points $y_{2:8}$, where $S_0(y_{2:k})$ has three connected components } 
\end{figure}
\end{center}
It results from the above that 
\begin{multline*}
k\mathbb{P}\left( \,\#\Phi_{\mathfrak{Q}_{\rho}}^{\eta }(\tau)=k\,\right) -\lambda _{d}(\mathfrak{Q}_{\rho})\cdot \mathbb{P}\left( \,\#\Phi_{\mathfrak{Q}_{\rho}}^{\eta \cup
\{0\}}(\tau)=k,g^{\eta \cup \{0\}}(0)>v_{\rho}(\tau)\,\right)   \\
=\frac{1}{(k-1)!}\sum_{m=1}^{k}\int_{\mathfrak{Q}_{\rho}}P_{x_{1}}[m]\mathrm{d}x_{1},
\end{multline*}
where, for any $x_{1}\in \mathbf{R}^{d}$, we write 
\begin{equation}
\label{eq:defPm}
P_{x_{1}}[m]:=\int_{(\mathfrak{Q}_{\rho}-x_{1})^{k-1}\cap
E_{m}}p_{x_{1}}(y_{2:k})\mathrm{d}y_{2:k}-\int_{\mathfrak{Q}_{\rho}^{k-1}\cap
E_{m}}p_{0}(y_{2:k})\mathrm{d}y_{2:k}.
\end{equation}
It is enough to show that for each $1\leq m\leq k$, we have  $\int_{\mathfrak{Q}_{\rho}}P_{x_{1}}[m]\mathrm{d}x_{1} = o\left(\lambda _{d}(\mathfrak{Q}_{\rho})\rho ^{-1}\right)$. To do it, we begin with $m=1$ which deals with
the case where there is exactly one connected component of size $k$ in $S_0(y_{2:k})$. Then we
extend our proof for $m\geq 2$ by dividing the set $S_0(y_{2:k})$ into its
connected components.

\paragraph{First case ($S_0(y_{2:k})$ has one connected component)} Assume that $m=1$ and let $x_{1}\in \mathbf{R}^{d}$ be fixed. We trivially obtain that 
\begin{align}
P_{x_{1}}[1]& =\int_{\left( (\mathfrak{Q}_{\rho}-x_{1})\cap \mathfrak{Q}_{\rho}\right) ^{k-1}\cap E_{1}}\left( p_{x_{1}}(y_{2:k})-p_{0}(y_{2:k})\right)
\mathrm{d}y_{2:k}  \notag  \label{eq:casemequal1} \\
& +\int_{(\mathfrak{Q}_{\rho}-x_{1})^{k-1}\cap E_{1}}p_{x_{1}}(y_{2:k})\mathbb{I}_{\exists j\leq k:y_{j}\in \mathfrak{Q}_{\rho}^{c}}\,\mathrm{d}y_{2:k} \\
& -\int_{\mathfrak{Q}_{\rho}^{k-1}\cap E_{1}}p_{0}(y_{2:k})\mathbb{I}_{\exists
j\leq k:y_{j}\in (\mathfrak{Q}_{\rho}-x_{1})^{c}}\,\mathrm{d}y_{2:k}.  \notag
\end{align}
We provide below a suitable upper bound for each term considered in the right-hand side of the above equation. 

\subparagraph{Upper bound for the first term in \eqref{eq:casemequal1}} 
Using the fact that for any events $A,B,C$, we have $|\PPP{A\cap B}-\PPP{A\cap C}|\leq \PPP{A\cap B^c\cap C}+\PPP{A\cap B\cap C^c}$, we obtain
\begin{multline}
\left\vert \int_{\left( (\mathfrak{Q}_{\rho}-x_{1})\cap \mathfrak{Q}_{\rho}\right)
^{k-1}\cap E_{1}}\left( p_{x_{1}}(y_{2:k})-p_{0}(y_{2:k})\right)
\mathrm{d}y_{2:k}\right\vert  \label{eq:bound1} \\
\begin{split}
& \leq \int_{\left( (\mathfrak{Q}_{\rho}-x_{1})\cap \mathfrak{Q}_{\rho}\right)
^{k-1}\cap E_{1}}\mathbb{P}\left( \,g^{\eta \cup
\{0,y_{2:k}\}}(0,y_{2:k})>v_{\rho}(\tau),M_{\mathfrak{Q}_{\rho}\setminus (\mathfrak{Q}_{\rho}-x_{1})}^{\eta \cup \{0,y_{2:k}\}}>v_{\rho}(\tau)\,\right) \mathrm{d}y_{2:k} \\
& +\int_{\left( (\mathfrak{Q}_{\rho}-x_{1})\cap \mathfrak{Q}_{\rho}\right)
^{k-1}\cap E_{1}}\mathbb{P}\left( \,g^{\eta \cup
\{0,y_{2:k}\}}(0,y_{2:k})>v_{\rho}(\tau),M_{(\mathfrak{Q}_{\rho}-x_{1})\setminus 
\mathfrak{Q}_{\rho}}^{\eta \cup \{0,y_{2:k}\}}>v_{\rho}(\tau)\,\right) \mathrm{d}y_{2:k}.
\end{split}
\end{multline}

To deal with the first term of the right-hand side of \eqref{eq:bound1}, we
introduce the event: 
\begin{equation*}
E^{\varnothing}(x_{1}):=\left\{ y_{2:k}\in \mathbf{R}^{(k-1)d}: 
S_0(y_{2:k}) \cap  \left(\mathfrak{Q}_{\rho}\setminus (\mathfrak{Q}_{\rho}-x_{1})\right) = \varnothing \right\}.
\end{equation*}
First, we assume that $y_{2:k}\in E^{\varnothing}(x_{1})\cap E_{1}$. From \eqref{eq:notalgebra}, we know that  $\left\{ g^{\eta \cup
\{0,y_{2:k}\}}(0,y_{2:k})>v_{\rho}(\tau)\right\}\in \Sigma_{\mathcal{A}}^{\eta\cup\{0,y_{2:k}\}}$ and 
$\{M_{\mathfrak{Q}_{\rho}\setminus (\mathfrak{Q}_{\rho}-x_{1})}^{\eta \cup \{0,y_{2:k}\}}>v_{\rho}(\tau)\}\in \Sigma_{\mathcal{B}}^{\eta\cup\{0,y_{2:k}\}}$, where \[\mathcal{A}:=\{\mathfrak{i}\in \mathcal{V}_{\rho}: \{0,y_{2:k}\}\cap \mathfrak{i}\neq \varnothing\} \quad \text{and}\quad  \mathcal{B}:=\{\mathfrak{j}\in \mathcal{V}_{\rho}: (\mathfrak{Q}_{\rho}\setminus (\mathfrak{Q}_{\rho}-x_{1}))\cap \mathfrak{j}\neq \varnothing \}.\] Since $y_{2:k}\in E^{\varnothing}(x_{1})$, we have $d(\mathcal{A}, \mathcal{B})>D$. It follows from Lemma \ref{Le:Arho}, \eqref{Le:Arho1} that,  conditional on  $\mathscr{A}_{\rho }$, the events $\left\{ g^{\eta \cup
\{0,y_{2:k}\}}(0,y_{2:k})>v_{\rho}(\tau)\right\} $ and $\left\{ M_{\mathfrak{Q}
_{\rho}\setminus (\mathfrak{Q}_{\rho}-x_{1})}^{\eta \cup \{0,y_{2:k}\}}>v_{\rho}(\tau)\right\} $ are independent.  This implies
that 
\begin{multline*}
\mathbb{P}\left( \,g^{\eta \cup \{0,y_{2:k}\}}(0,y_{2:k})>v_{\rho}(\tau),M_{\mathfrak{Q}_{\rho}\setminus (\mathfrak{Q}_{\rho}-x_{1})}^{\eta \cup
\{0,y_{2:k}\}}>v_\rho(\tau)\,\right) \mathbb{I}_{y_{2:k}\in E^{\varnothing}(x_{1})}\,  \\
\leq \PPP{\mathscr{A}_{\rho}}^{-1}\cdot \mathbb{P}\left( \,g^{\eta \cup \{0,y_{2:k}\}}(0,y_{2:k})>v_{\rho}(\tau)\,\right) \cdot  \mathbb{P}\left( \,M_{\mathfrak{Q}_{\rho}\setminus (\mathfrak{
Q}_{\rho}-x_{1})}^{\eta \cup \{0,y_{2:k}\}}>v_\rho(\tau)\,\right) +\mathbb{P}\left( \,(\mathscr{A}_{\rho })^{c}\,\right),
\end{multline*}
where we have bounded the indicator function $\mathbb{I}_{y_{2:k}\in E^{\varnothing}(x_{1})}$ by 1. According to Lemma \ref{Le:Arho}, \eqref{Le:Arho2}, we know that  $\PPP{\mathscr{A}_{\rho}} \geq 1 - c\cdot \rho^{-\alpha}$ for any $\alpha>0$. Moreover, since $g$ satisfies Condition (C), we  have $\mathbb{P}\left( \,g^{\eta \cup \{0,y_{2:k}\}}(0,y_{2:k})>v_{\rho}(\tau)\,\right) \leq c\cdot \rho^{-1}$.  Besides, it results from the Slivnyak-Mecke formula that 
\begin{multline*}
\mathbb{P}\left( \,M_{\mathfrak{Q}_{\rho}\setminus (\mathfrak{Q}
_{\rho}-x_{1})}^{\eta \cup \{0,y_{2:k}\}}>v_\rho(\tau)\,\right)\\
\begin{split}
 & = \PPP{\exists z\in (\mathfrak{Q}_{\rho}\setminus (\mathfrak{Q}
_{\rho}-x_{1}))\cap (\eta\cup\{0,y_{2:k}\}):  g^{\eta\cup\{0,y_{2:k}\}}(z)>v_\rho(\tau)}\\
& \leq \mathbb{E}\left[ \,\sum_{z\in \eta \cap \mathfrak{Q}_{\rho}}\mathbb{
I}_{g^{\eta \cup \{0,y_{2:k}\}}(z)>v_\rho(\tau)}\,\,\right] + \sum_{j=2}^{k}\mathbb{P}\left( \,g^{\eta \cup \{0,y_{2:k}\}}(y_{j})>v_{\rho}(\tau)\,\right) \\
& =  \int_{\mathfrak{Q}_{\rho}}\mathbb{P}\left( \,g^{\eta
\cup \{0,y_{2:k},z\}}(z)>v_\rho(\tau)\,\right) \mathrm{d}z + \sum_{j=2}^{k}\mathbb{P}\left( \,g^{\eta \cup \{0,y_{2:k}\}}(y_{j})>v_{\rho}(\tau)\,\right)    \\
& \leq c\cdot \lambda _{d}(\mathfrak{Q}_{\rho})\rho ^{-1},
\end{split}
\end{multline*}
where the last line is also a consequence of Condition (C). This implies that 
\begin{equation}
\label{eq:majE11}
\mathbb{P}\left( \,g^{\eta \cup \{0,y_{2:k}\}}(0,y_{2:k})>v_\rho(\tau),M_{\mathfrak{Q}_{\rho}\setminus (\mathfrak{Q}_{\rho}-x_{1})}^{\eta \cup
\{0,y_{2:k}\}}>v_\rho(\tau)\,\right) \mathbb{I}_{y_{2:k}\in E^{\varnothing}(x_{1})\cap E_{1}}\,\leq c\cdot \lambda _{d}(\mathfrak{Q}_{\rho})\rho ^{-2}.
\end{equation}
Secondly, we assume that $y_{2:k}\in (E^{\varnothing}(x_{1}))^{c}\cap E_{1}$. In
particular, we have 
\begin{equation*}
\delta (0,\mathfrak{Q}_{\rho}\setminus (\mathfrak{Q}_{\rho}-x_{1}))\leq \text{diam}
(S_0(y_{2:k}))\leq c\cdot \lambda _{d}(\mathfrak{c}_{\rho})^{1/d}.
\end{equation*}
Since $g$ satisfies Condition (C), this implies that 
\begin{multline}
\label{eq:majE12}
\mathbb{P}\left( \,g^{\eta \cup \{0,y_{2:k}\}}(0,y_{2:k})>v_\rho(\tau),M_{\mathfrak{Q}_{\rho}\setminus (\mathfrak{Q}_{\rho}-x_{1})}^{\eta \cup
\{0,y_{2:k}\}}>v_\rho(\tau)\,\right) \mathbb{I}_{y_{2:k}\in (E^{\varnothing
}(x_{1}))^{c}\cap E_{1}}\, \\
\leq c\cdot \rho ^{-1}\mathbb{I}_{\delta (0,\mathfrak{Q}_{\rho}\setminus (\mathfrak{Q}_{\rho}-x_1))\leq c\cdot \lambda _{d}(\mathfrak{c}_{\rho})^{1/d}}\,.
\end{multline}
Integrating over $y_{2:k}\in \left( (\mathfrak{Q}_{\rho}-x_{1})\cap \mathfrak{Q}_{\rho}\right) ^{k-1}\cap
E_{1}$, it follows from \eqref{eq:majE11} and \eqref{eq:majE12} that 
\begin{multline}
\int_{\left( (\mathfrak{Q}_{\rho}-x_{1})\cap \mathfrak{Q}_{\rho}\right) ^{k-1}\cap
E_{1}}\mathbb{P}\left( \,g^{\eta \cup \{0,y_{2:k}\}}(0,y_{2:k})>v_\rho(\tau),M_{\mathfrak{Q}_{\rho}\setminus (\mathfrak{Q}_{\rho}-x_{1})}^{\eta \cup
\{0,y_{2:k}\}}>v_\rho(\tau)\,\right) \mathrm{d}y_{2:k}  \label{eq:majcase11} \\
\begin{split}
& \leq c\cdot \lambda _{(k-1)d}\left( 
E_{1}\right) \cdot \left( \lambda _{d}(\mathfrak{Q}_{\rho})\rho ^{-2}+\rho ^{-1}\mathbb{I}_{\delta (0,\mathfrak{Q}_{\rho}\setminus (\mathfrak{Q}_{\rho}-x_{1}))\leq c\cdot \lambda _{d}(\mathfrak{c}_{\rho})^{1/d}}\,\right) \\
& \leq c\cdot \lambda _{d}(\mathfrak{c}_{\rho})^{k-1}\cdot \left( \lambda _{d}(
\mathfrak{Q}_{\rho})\cdot \rho ^{-2}+\rho ^{-1}\ind{\delta (0,\mathfrak{Q}_{\rho}\setminus (\mathfrak{Q}_{\rho}-x_{1}))\leq c\cdot \lambda _{d}(\mathfrak{c}_{\rho})^{1/d}}\,\right),
\end{split}
\end{multline} where the last line comes from the fact that  $\lambda_{(k-1)d}(E_1)\leq c\cdot \lambda_d(\mathfrak{c}_{\rho})^{k-1}$. 

Proceeding exactly along the same lines as above, by considering here the event:
\begin{equation*}
F^{\varnothing}(x_{1}):=\left\{ y_{2:k}\in \mathbf{R}^{(k-1)d}: 
S_0(y_{2:k}) \cap  \left( (\mathfrak{Q}_{\rho}-x_{1})\setminus\mathfrak{Q}_{\rho}\right) = \varnothing \right\},
\end{equation*}
 we can show that the second term of the right-hand side of \eqref{eq:bound1} can be bounded in a similar way, i.e. 
\begin{multline*}
\int_{\left( (\mathfrak{Q}_{\rho}-x_{1})\cap \mathfrak{Q}_{\rho}\right) ^{k-1}\cap
E_{1}}\mathbb{P}\left( \,g^{\eta \cup \{0,y_{2:k}\}}(0,y_{2:k})>v_{\rho}(\tau),M_{(\mathfrak{Q}_{\rho}-x_{1})\setminus \mathfrak{Q}_{\rho}}^{\eta \cup
\{0,y_{2:k}\}}>v_\rho(\tau)\,\right) \mathrm{d}y_{2:k}   \\
\leq c\cdot \lambda _{d}(\mathfrak{c}_{\rho})^{k-1}\cdot \left( \lambda _{d}(
\mathfrak{Q}_{\rho})\cdot \rho ^{-2}+\rho ^{-1}\ind{\delta (0,(\mathfrak{Q}_{\rho}-x_{1})\setminus \mathfrak{Q}_{\rho}))\leq c\cdot \lambda _{d}(\mathfrak{c}_{\rho})^{1/d}}\,\right) .
\end{multline*}

This together with \eqref{eq:bound1} and \eqref{eq:majcase11} implies that 
\begin{multline*}\int_{\left( (\mathfrak{Q}_{\rho}-x_{1})\cap \mathfrak{Q}_{\rho}\right) ^{k-1}\cap E_{1}}\left( p_{x_{1}}(y_{2:k})-p_{0}(y_{2:k})\right)
\mathrm{d}y_{2:k}\leq  c\cdot \lambda _{d}(\mathfrak{c}_{\rho})^{k-1}\cdot  \lambda _{d}(
\mathfrak{Q}_{\rho})\cdot \rho ^{-2}\\
+ c\cdot \lambda _{d}(\mathfrak{c}_{\rho})^{k-1}\cdot\rho^{-1}\cdot \left( \ind{\delta (0,\mathfrak{Q}_{\rho}\setminus (\mathfrak{Q}_{\rho}-x_{1}))\leq c\cdot \lambda _{d}(\mathfrak{c}_{\rho})^{1/d}} + \ind{\delta (0,(\mathfrak{Q}_{\rho}-x_{1})\setminus \mathfrak{Q}_{\rho}))\leq c\cdot \lambda _{d}(\mathfrak{c}_{\rho})^{1/d}}\right).
 \end{multline*}
This deals with the first term of the right-hand side in  \eqref{eq:casemequal1}.

\subparagraph{Upper bound for the second term in \eqref{eq:casemequal1}}
Trivially, this term equals 0: this comes from the fact that for $\rho$ large enough, we have 
\begin{equation*}
E_1\cap \{y_{2:k}\in \mathbf{R}^{(k-1)d}: \exists j\leq k \text { s.t. }
y_j\in \mathfrak{Q}_{\rho}^c \} = \varnothing
\end{equation*}
since $\lambda_d(\mathfrak{c}_{\rho}) = o\left( \lambda_d(\mathfrak{Q}_{\rho}) \right)$. 

\subparagraph{Upper bound for the third term in \eqref{eq:casemequal1}} To deal
with this term, we notice that if 
$E_{1}\cap \left\{ y_{2:k}\in \mathfrak{Q}_{\rho}^{k-1}:\exists j\leq k\text{
s.t. }y_{j}\in (\mathfrak{Q}_{\rho}-x_{1})^{c}\right\} \neq \varnothing$, then 
$\delta \left( 0,\mathfrak{Q}_{\rho}\setminus (\mathfrak{Q}_{\rho}-x_{1})\right) \leq c\cdot \lambda _{d}(\mathfrak{c}_{\rho})^{1/d}$ because $\text{diam}(S_0(y_{2:k}))\leq c\cdot \lambda _{d}(\mathfrak{c}_{\rho})^{1/d}$. Besides, since $p_{0}(y_{2:k})\leq c\cdot \rho ^{-1}$ according to Condition (C), we obtain by integrating over $y_{2:k}\in \mathfrak{Q}_{\rho}^{k-1}\cap E_{1}$ that 
\begin{equation*}
\int_{\mathfrak{Q}_{\rho}^{k-1}\cap E_{1}}p_{0}(y_{2:k})
\mathbb{I}_{\exists j\leq k:y_{j}\in (\mathfrak{Q}_{\rho}-x_{1})^{c}}
\,\mathrm{d}y_{2:k}\leq c\cdot \lambda _{d}(\mathfrak{c}_{\rho})^{k-1}\cdot \rho^{-1}\cdot \ind{\delta \left( 0,\mathfrak{Q}_{\rho}\setminus (\mathfrak{Q}_{\rho}-x_{1})\right) \leq c\cdot \lambda _{d}(\mathfrak{c}_{\rho})^{1/d}}. \label{eq:casemequal1case3}
\end{equation*}
This deals with the third term of the right-hand side in  \eqref{eq:casemequal1}.

By considering the three upper bounds discussed above and by integrating over $x_1\in \mathfrak{Q}_{\rho}$, we get 
\begin{multline*}
\left|\int_{\mathfrak{Q}_{\rho}} P_{x_1}[1]\mathrm{d}x_1\right| \leq c\cdot \lambda _{d}(\mathfrak{c}_{\rho})^{k-1}\cdot  \lambda _{d}(
\mathfrak{Q}_{\rho})^2\cdot \rho ^{-2}\\
\begin{split}
& + c\cdot \lambda _{d}(\mathfrak{c}_{\rho})^{k-1}\cdot \rho^{-1}\cdot \lambda _{d}\left( \left\{x_{1}\in \mathfrak{Q}_{\rho}:\delta (0,\mathfrak{Q}_{\rho}\setminus (\mathfrak{Q}_{\rho}-x_{1})\leq c\cdot \lambda_d(\mathfrak{c}_{\rho})^{1/d}\right\}\right)\\
&  + c\cdot \lambda _{d}(\mathfrak{c}_{\rho})^{k-1}\cdot \rho^{-1}\cdot \lambda _{d}\left( \left\{x_{1}\in \mathfrak{Q}_{\rho}:\delta (0, (\mathfrak{Q}_{\rho}-x_{1})\setminus \mathfrak{Q}_{\rho}\leq c\cdot \lambda_d(\mathfrak{c}_{\rho})^{1/d}\right\}\right).
\end{split}
\end{multline*}
We can easily prove that if $x_1\in  \mathfrak{Q}_{\rho}$ is such that $\delta (0,\mathfrak{Q}_{\rho}\setminus (\mathfrak{Q}_{\rho}-x_{1})\leq c\cdot \lambda _{d}(\mathfrak{c}_{\rho})^{1/d}$, then $x_1\in \mathfrak{Q}_{\rho}\setminus (\mathfrak{Q}_{\rho} \ominus c^{1/d}\cdot \mathfrak{c}_{\rho})$. Hence 
\begin{equation*}
\begin{split}
\lambda _{d}\left( \left\{x_{1}\in \mathfrak{Q}_{\rho}:\delta (0,\mathfrak{Q}_{\rho}\setminus (\mathfrak{Q}_{\rho}-x_{1})\leq c\cdot \lambda _{d}(\mathfrak{c}_{\rho})^{1/d}\right\}\right)
& \leq \lambda_d\left(\mathfrak{Q}_{\rho}\setminus (\mathfrak{Q}_{\rho} \ominus c^{1/d}\cdot \mathfrak{c}_{\rho})\right)\\
 & = O\left( \lambda _{d}(\mathfrak{Q}_{\rho})^{(d-1)/d}\cdot \lambda_d(\mathfrak{c}_{\rho})^{1/d}\right).  
 \end{split}
\end{equation*}
Moreover, for $\rho$ large enough, we have
\begin{equation*}
\lambda _{d}\left( \{x_{1}\in \mathfrak{Q}_{\rho}:\delta (0,    (\mathfrak{Q}_{\rho}-x_{1})\setminus \mathfrak{Q}_{\rho}\leq c\cdot \lambda _{d}(\mathfrak{c}
_{0})^{1/d}\}\right) =0
\end{equation*}
since $\lambda_d(\mathfrak{c}_{\rho}) = o\left( \lambda_d(\mathfrak{Q}_{\rho})\right)$. From \eqref{def:qrho}, we deduce that 
\[\left| \int_{\mathfrak{Q}_{\rho}}P_{x_{1}}[1]\mathrm{d}x_{1}\right| \leq c\cdot \left(\lambda_d(\mathfrak{c}_\rho)^{k-1+1/d}\cdot \lambda_d(\mathfrak{Q}_\rho)^{-1/d} \right)\cdot \left(\lambda _{d}(\mathfrak{Q}_{\rho})\cdot\rho ^{-1}\right) =o\left( \lambda _{d}(\mathfrak{Q}_{\rho})\cdot\rho ^{-1}\right).\] This concludes the proof for the case $m=1$.

\paragraph{Second case ($S_0(y_{2:k})$ has $m$ connected component with $m\geq 2$ )} Assume that $m\geq 2$ and $y_{2:k}\in E_m$.

First, we provide below a uniform upper bound for $p_{x}(y_{2:k})$, with $
x\in \mathfrak{Q}_{\rho}$. Since $y_{2:k}\in E_{m}$, we can divide $S_0(y_{2:k})$
into its $m$ connected components, say $C_{1}(y_{2:k}),\ldots
,C_{m}(y_{2:k}) $. For each $1\leq l\leq m$, let $J_{l}\subset \{1,\ldots
,k\}$ be the set of indices $j$ such that $C_{l}(y_{2:k})=\bigcup_{j\in
J_{l}}(y_{j}+\mathfrak{c}_{\rho})$, with  $y_1:=0$. In particular, we have 
\begin{equation*}
p_{x}(y_{2:k}) \leq \PPP{g^{\eta\cup\{0, y_{2:k}\}}(0,y_{2:k})>v_\rho(\tau)}  =\mathbb{P}\left( \,\bigcap_{l=1}^{m}\left\{ g^{\eta \cup
\{0,y_{2:k}\}}(y_{J_{l}}) >v_\rho(\tau)\right\}\,\right) ,
\end{equation*}
where we recall that $g^{\eta \cup \{0,y_{2:k}\}}(y_{J_{l}})>v_\rho(\tau)$
means that $g^{\eta \cup \{0,y_{2:k}\}}(y)>v_\rho(\tau)$ for any $y\in J_{l}$. In the same spirit as in the case where $S_0(y_{2:k})$ has one connected component, we deduce from  Lemma \ref{Le:Arho}, \eqref{Le:Arho1} that,  conditional on the event $\mathscr{A}_{\rho }$, the events $\{g^{\eta \cup
\{0,y_{2:k}\}}(y_{J_{l}})>v_\rho(\tau)\}$, $1\leq l\leq m$, are independent. This gives 
\begin{equation*}
p_{x}(y_{2:k})\leq \PPP{\mathscr{A}_{\rho}}^{-(m-1)}\cdot \prod_{l=1}^{m}\mathbb{P}\left( \,g^{\eta \cup
\{0,y_{2:k}\}}(y_{J_{l}})>v_\rho(\tau)\,\right) + \PPP{\mathscr{A}_{\rho }^{c}} .
\end{equation*}
Since $g$ satisfies Condition (C), it follows from Lemma \ref{Le:Arho}, \eqref{Le:Arho2} that there
exists a constant $c>0$ such that, for any $x\in \mathfrak{Q}_{\rho}$ and for
any $y_{2:k}\in E_{m}$, we have $p_{x}(y_{2:k})\leq c\cdot \rho ^{-m}$.

Now, we are able to provide an upper bound for $\left\vert \int_{\mathfrak{Q}_{\rho}}P_{x_{1}}[m]\mathrm{d}x_{1}\right\vert $. Indeed, integrating over $y_{2:k}$ in the right-hand side of \eqref{eq:defPm}, we have 
\begin{equation*}
|P_{x_{1}}[m]|\leq c\cdot \rho ^{-m}\cdot \sup_{x\in \mathfrak{Q}_{\rho}}\lambda _{(k-1)d}((\mathfrak{Q}_{\rho}-x)^{k-1}\cap E_{m})\leq c\cdot \rho
^{-m}\cdot \lambda _{d}(\mathfrak{Q}_{\rho})^{m-1}\cdot \lambda _{d}(\mathfrak{c}_{\rho})^{k-m}.
\end{equation*}
Integrating over $x_{1}\in \mathfrak{Q}_{\rho}$, we get 
\begin{equation*}
\left\vert \int_{\mathfrak{Q}_{\rho}}P_{x_{1}}[m]\mathrm{d}x_{1}\right\vert \leq c\cdot 
\left(\rho^{-1}\cdot \lambda_d(\mathfrak{Q}_\rho)\cdot \lambda_d(\mathfrak{c}_\rho)^{(k-m)/(m-1)}    \right)^{m-1}\cdot \left(\lambda_d(\mathfrak{Q}_{\rho})\cdot \rho^{-1}\right)  = o\left(\lambda_d(\mathfrak{Q}_{\rho})\cdot \rho^{-1} \right). 
\end{equation*}
This concludes the proof of Proposition \ref{Prop:identification} for any $ k\geq 2$. The case $k=1$ is much more simple than the case $k\geq 2$ and can be dealt by following the same lines as above and by noting that 
\begin{equation*}
\mathbb{P}\left( \,\#\Phi_{\mathfrak{Q}_{\rho}}^{\eta }(\tau )=1\,\right) -\lambda
_{d}(\mathfrak{Q}_{\rho})\cdot \mathbb{P}\left( \,\#\Phi_{\mathfrak{Q}_{\rho}}^{\eta
\cup \{0\}}(\tau )=1,g^{\eta \cup \{0\}}(0)>v_{\rho }(\tau )\,\right) =\int_{\mathfrak{Q}_{\rho}}(p_{x_{1}}-p_{0})\mathrm{d}x_1,
\end{equation*}
where, for any $x\in \mathfrak{Q}_{\rho}$, we write $p_{x}:=\mathbb{P}\left(
\,g^{\eta \cup \{0\}}(0)>v_\rho(\tau),M^{\eta\cup\{0\}}_{(\mathfrak{Q}_{\rho}-x)\setminus
\{0\}}\leq v_\rho(\tau)\,\right) $.
\end{prooft}

\subsection{Our main theorem}

Let $g$ be a geometric characteristic such that \eqref{Threshcond} holds for
some $\tau_0 >0$. According to Leadbetter \cite{L2}, we say that the extremal index $\theta \in [0,1]$ of the Poisson-Voronoi tessellation  exists if $\lim_{\rho\rightarrow \infty }\mathbb{P}(\#
\Phi^\eta _{\mathbf{W}_\rho}(\tau_0)=0)=e^{-\theta \tau_0 }$. We are now prepared to state our
main theorem on the weak convergence of the point process $\Phi^\eta _{\mathbf{W}_\rho}(\tau)$ for each $\tau>0$.   

\begin{Th}
\label{Th:theta} Let $g$ be a geometric characteristic satisfying Condition (C). Assume that there
exist $\tau _{0}>0$ such that \eqref{Threshcond} holds and $(a_{k})_{k\geq 1}$ such 
that $\pi _{k,\mathfrak{Q}_{\rho}}(\tau_0)\leq a_{k}$ for any $k\geq 1$ and any $\rho >0$, with $\sum_{k=1}^{\infty }a_{k}<\infty$. 

\begin{enumerate}[(i)]
\item The following assertions are equivalent:
\begin{enumerate}[(A)]
\item there exists $\theta \in (0,1]$ such that $\lim_{\rho\rightarrow \infty }\mathbb{P}(\#
\Phi^\eta _{\mathbf{W}_\rho}(\tau_0)=0)=e^{-\theta \tau_0 }$ and the following limit
exist $p_{k}:=\lim_{\rho \rightarrow \infty }p_{k,\mathfrak{Q}_{\rho}}(\tau_0)$ for any 
$k\geq 1$;\label{Th:Assert1}

\item for any $\tau >0$, the point process $\Phi^\eta _{\mathbf{W}_\rho}(\tau)$ converges to a homogeneous compound Poisson point process in $W:=\left[ -1/2,1/2\right] ^{d}$ with intensity $\nu(\tau)
>0$ and cluster size distributions $\pi _{k}:=\lim_{\rho \rightarrow \infty
}\pi _{k,\mathfrak{Q}_{\rho}}(\tau_0)$, with $k\geq 1$.\label{Th:Assert2}
\end{enumerate}
\item If one of the above assertions holds, we have $p_{k}=k\theta \pi _{k}$ for any $k\geq 1$ and $\theta =\sum_{k=1}^{\infty }k^{-1}p_{k}$.
\end{enumerate}

\end{Th}

Our theorem provides a new characterization of the extremal index. Indeed, this index was previously interpreted as the reciprocal of the mean of the
cluster size distribution $\pi$. Now, it can be viewed as the mean of the
reciprocal of the Palm version of the cluster size. Besides, our new characterization: $\theta =\sum_{k=1}^{\infty }k^{-1}p_{k}$ will be extensively used in Section \ref{sec:examples} to estimate the extremal indices for various geometric characteristics. 

To prove Theorem \ref{Th:theta}, we associate with the point process $\Phi^\eta _{\mathbf{W}_\rho}(\tau)$ its Laplace transform $\mathbb{L}_{\rho }$ defined as follows: for any
continuous function $f:W\rightarrow \mathbf{R}_{+}$, we have  
\begin{equation*}
\begin{tabular}{l}
$\mathbb{L}_{\rho}(f):=\mathbb{E}\left[ \,\exp \left( -\sum_{y\in \Phi^\eta _{\mathbf{W}_\rho}(\tau)}f(y)\right) \,\right]$. 
\end{tabular}
\end{equation*}
It is well-known that the weak convergence of $\Phi^\eta _{\mathbf{W}_\rho}(\tau)$ is
equivalent to the convergence of its Laplace transform for any positive and continuous
function $f$. For a sequence of real random variables, the weak convergence
of the point process of exceedances has been investigated in \cite{HHL} and
generalized to random fields on $\mathbf{N}_{+}^{d}$ in \cite{FP}. We use below 
the same type of approach. However, we have to take into account specific features
of random tessellations in $\mathbf{R}^{d}$.

The first step consists in showing that exceedances over disjoint sub-cubes
behave asymptotically as if they were independent. To do it, we divide 
$\mathbf{W}_\rho$ into $m_{\rho}^{d}$ disjoint sub-cubes $\mathfrak{B}[l]$, $l=1,\ldots
,m_{\rho}^{d}$, with the same volume as $\mathfrak{Q}_{\rho}$, where $m_{\rho}$ is defined in \eqref{def:nm}. 

\begin{Le}
\label{Le:independencePP}

\begin{enumerate}[(i)]
\item \label{independencePP1}For any measurable function $f:W\rightarrow 
\mathbf{R}_{+}$, we have 
\begin{equation*}
\mathbb{L}_{\rho}(f)-\prod_{l=1}^{m_{\rho}^{d}}\mathbb{E}\left[ \,\exp \left( -\sum_{y\in
\Phi^\eta _{\mathbf{W}_\rho}(\tau)\cap \rho ^{-1/d}\mathfrak{B}[l]}f(y)\right) \,\right] 
\underset{\rho \rightarrow \infty }{\longrightarrow }0.
\end{equation*}

\item \label{independencePP2}Moreover, we have
\begin{equation*}
\mathbb{P}\left( M_{\mathbf{W}_\rho}^{\eta }\leq v_{\rho }(\tau )\right)
-\prod_{l=1}^{m_{\rho}^{d}}\mathbb{P}\left( M_{\mathfrak{B}[l]}^{\eta }\,\leq
v_{\rho }(\tau )\right) \underset{\rho \rightarrow \infty }{\longrightarrow }0.
\end{equation*}
\end{enumerate}
\end{Le}


\begin{prooft}
We begin with the first assertion. For any $l\leq m_{\rho}^d$, we write  $\mathfrak{B}^{\circ}[l]:=\mathfrak{B}[l]\ominus \mathfrak{c}_{\rho}$.  Let
\begin{equation*}
L_{\rho ,l}(f)\,=\exp \left( -\sum_{y\in \Phi^\eta_{\mathbf{W}_\rho}(\tau)\cap \rho ^{-1/d}\mathfrak{B}[l]}f(y)\right)
\quad \text{and}\quad 
L_{\rho ,l}^{\circ }(f)=\exp \left( -\sum_{y\in \Phi^\eta_{\mathbf{W}_\rho}(\tau)\cap \rho ^{-1/d}\mathfrak{B}^{\circ}[l]}f(y)\,\right).
\end{equation*}
We write
\begin{equation*}
\mathbb{L}_{\rho}(f)-\prod_{l=1}^{m_{\rho}^{d}}\mathbb{E}\left[ L_{\rho
,l}(f)\right] =\Delta L_{\rho ,1}(f)+\Delta L_{\rho ,2}(f)+\Delta L_{\rho
,3}(f)+\Delta L_{\rho ,4}(f),
\end{equation*}
where
\begin{equation*}
\begin{tabular}{l}
$\Delta L_{\rho ,1}(f)=\mathbb{E}\left[ \prod_{l=1}^{m_{\rho}^{d}}L_{\rho ,l}(f)\,\right] -\mathbb{E}\left[ \prod_{l=1}^{m_{\rho}^{d}}L_{\rho ,l}^{\circ }(f)\,\right] ,$ \\ 
$\Delta L_{\rho ,2}(f)=\mathbb{E}\left[ \prod_{l=1}^{m_{\rho}^{d}}L_{\rho ,l}^{\circ }(f)\,\right] -\prod_{l=1}^{m_{\rho}^{d}}\mathbb{E}\left[ \,L_{\rho ,l}^{\circ }(f)\,\right] ,$ \\ 
$\Delta L_{\rho ,3}(f)=\prod_{l=1}^{m_{\rho}^{d}}\mathbb{E}\left[ \,L_{\rho ,l}^{\circ }(f)\,\right] -\prod_{l=1}^{m_{\rho}^{d}}\mathbb{E}
\left[ L_{\rho ,l}(f)\,\right] ,$ \\ 
$\Delta L_{\rho ,4}(f)=\mathbb{E}\left[ \,\exp \left( -\sum_{y\in \Phi^\eta_{\mathbf{W}_\rho}(\tau)}f(y)\right) \,\right] -\mathbb{E}\left[ \prod_{l=1}^{m_{\rho}^{d}}
L_{\rho ,l}(f)\,\right] .$\end{tabular}
\end{equation*}

We prove below that each term converges to $0$. For
the third term, using the fact that $|\prod x_{i}-\prod y_{i}|\leq \sum
|x_{i}-y_{i}|$ for $0\leq x_{i},y_{i}\leq 1$ and the fact that $|\exp (-x)-\exp (-y)|\leq
|x-y|$ for all $x,y\geq 0$, we get 
\begin{align*}
|\Delta L_{\rho ,3}(f)|& \leq m_{\rho}^{d}\cdot \sup_{l\leq m_{\rho}^{d}}\mathbb{E}\left[
\,\sum_{y\in \Phi^\eta _{\mathbf{W}_\rho}(\tau)\cap \rho ^{-1/d}(\mathfrak{B}[l]\backslash 
\mathfrak{B}^{\circ}[l])}f(y)\right]   \notag \\
& \leq m_{\rho}^{d}\cdot \sup_{l\leq m_{\rho}^{d}}\mathbb{E}\left[ \,\sum_{x\in \eta \cap
\left( \mathfrak{B}[l]\backslash \mathfrak{B}^{\circ}[l]\right) }f(\rho
^{-1/d}x)\mathbb{I}_{g^\eta(x)>v_{\rho }(\tau)}\right]  \notag\\
& = c\cdot m_{\rho}^{d}\cdot \sup_{l\leq m_{\rho}^{d}}\int_{\mathfrak{B}[l]\backslash 
\mathfrak{B}^{\circ}[l]}f(\rho ^{-1/d}x)\mathbb{P}\left( \,g^{\eta \cup
\{x\}}(x)>v_{\rho }(\tau)\,\right) \mathrm{d}x \notag \\
& \leq c\cdot m_{\rho}^{d}\cdot \lambda _{d}\left( \mathfrak{B}[l]\backslash 
\mathfrak{B}^{\circ}[l]\right) \cdot \mathbb{P}\left( \,g(\mathcal{C})>v_{\rho }(\tau)\right) \,,
\end{align*}
where the third line comes from the Slivnyak-Mecke formula and where the
fourth line comes from \eqref{def:typicalcell} and the fact that $f$ is
bounded because it is continuous on the compact set $W$. Since $m_{\rho}^{d}\underset{\rho \rightarrow \infty }{\sim }\rho \cdot q_{\rho }^{-1}\cdot (\log \rho
)^{-(1+\varepsilon) }$ and
\begin{equation*}
\lambda _{d}\left( \mathfrak{B}[l]\backslash \mathfrak{B}^{\circ}[l]\right) \leq c\cdot q_{\rho }^{(d-1)/d}\cdot (\log \rho )^{(1+\varepsilon
)},
\end{equation*}
we deduce that
\begin{equation*}
|\Delta L_{\rho ,3}(f)| = O\left(q_{\rho }^{-1/d}\right).
\end{equation*}
In the same spirit as above, we prove that $\Delta L_{\rho ,1}(f)$ and $\Delta L_{\rho ,4}(f)$ converges to $0$.

For $\Delta L_{\rho ,2}(f)$, we notice that conditional on $\mathscr{A}_{\rho }$, the
random variables considered in the expectations are independent. Then, we
have
\begin{equation*}
\mathbb{E}\left[ \prod_{l=1}^{m_{\rho}^{d}}L_{\rho ,l}^{\circ
}(f)\,\right] =\prod_{l=1}^{m_{\rho}^{d}}\mathbb{E}\left[ \left. \,L_{\rho ,l}^{\circ }(f)\right\vert \,\mathscr{A}_{\rho }\right] +\mathbb{P}\left( \,\mathscr{A}_{\rho }^{c}\right) \left( \mathbb{E}\left[ \left.
\prod_{l=1}^{m_{\rho}^{d}}L_{\rho ,l}^{\circ }(f)\right\vert
\,\mathscr{A}_{\rho }^{c}\right] -\prod_{l=1}^{m_{\rho}^{d}}\mathbb{E}\left[ \left. \,
L_{\rho ,l}^{\circ }(f)\right\vert \,\mathscr{A}_{\rho }\right]
\right) .
\end{equation*}
Moreover
\begin{eqnarray*}
\prod_{l=1}^{m_{\rho}^{d}}\mathbb{E}\left[ \,L_{\rho ,l}^{\circ
}(f)\,\right]  &=&\prod_{l=1}^{m_{\rho}^{d}}\left( \mathbb{E}\left[ \left. \,
L_{\rho ,l}^{\circ }(f)\right\vert \,\mathscr{A}_{\rho }\right] +
\mathbb{P}\left( \,\mathscr{A}_{\rho }^{c}\right) \left( \mathbb{E}\left[ \left. \,
L_{\rho ,l}^{\circ }(f)\right\vert \,\mathscr{A}_{\rho }^{c}\right] -
\mathbb{E}\left[ \left. \,L_{\rho ,l}^{\circ
}(f)\right\vert \,\mathscr{A}_{\rho }\right] \right) \right)  \\
&:=&\prod_{l=1}^{m_{\rho}^{d}}\mathbb{E}\left[ \left. \,L_{\rho
,l}^{\circ }(f)\right\vert \,\mathscr{A}_{\rho }\right] +m_{\rho}^{d}\mathbb{P}\left(
\,\mathscr{A}_{\rho }^{c}\right) H_{\rho }(f).
\end{eqnarray*}
The term $H_{\rho }(f)$ appearing in the above equation is such that $|H_{\rho }(f)|\leq c$: this is a consequence of Lemma \ref{Le:Arho}, \eqref{Le:Arho2} and the fact that $ 0\leq \EEE{\left. L_{\rho ,l}^{\circ }(f)\right\vert \,\mathscr{A}_{\rho }^{c}}\leq 1$ and $0\leq \EEE{ \left. \,L_{\rho ,l}^{\circ
}(f)\right\vert \,\mathscr{A}_{\rho }}\leq 1$. By applying again Lemma \ref{Le:Arho}, \eqref{Le:Arho2}, it follows that $\Delta L_{\rho ,2}(f)$ converges
to $0$.  We proceed in a similar way for the proof of the second assertion.
\end{prooft}


We now adapt two theorems due to Leadbetter, Lindgren and Rootz\'{e}n in our context. The following result is an adaptation of Theorem 4.2 in \cite{HHL} (resp. Proposition 4.2 in \cite{FP}) and gives sufficient conditions to derive the convergence of $\Phi^\eta _{\mathbf{W}_\rho}(\tau)$ to a homogeneous compound Poisson point process.

\begin{Prop}
\label{Prop:HHL1} Assume that $\mathbb{P}\left( \,\# \Phi^\eta _{\mathbf{W}_\rho}(\tau_0)=0\right) \underset{\rho \rightarrow \infty }{\longrightarrow }e^{-\nu }$
for some $\tau _{0}>0$ and $\nu >0$. If $\left( \pi _{k,\mathfrak{Q}_{\rho}}\right) _{k\geq 1}$ converges to a probability distribution $\pi $ on $\mathbf{N}_{+}$, then $\Phi^\eta _{\mathbf{W}_\rho}(\tau _{0})$ converges in distribution to a homogeneous compound Poisson point process with intensity $\nu $ and
limiting cluster size distribution $\pi $.
\end{Prop}

The following result adapted from Theorem  5.1 in \cite{HHL} (resp. Proposition
4.3. in \cite{FP}) shows that if $\Phi^\eta _{\mathbf{W}_\rho}(\tau_0)$ has a limit for some 
$\tau_0>0$, it has a limit for all $\tau >0$.

\begin{Prop}
\label{Prop:HHL2} Assume that $\Phi^\eta _{\mathbf{W}_\rho}(\tau_0)$ converges to a homogeneous compound
Poisson point process in $W$ with intensity $\nu >0$ and cluster size
distribution $\pi $, for some $\tau_0>0$. Then $\Phi^\eta _{\mathbf{W}_\rho}(\tau)$
converges to a homogeneous compound Poisson point process with intensity $\nu \cdot \tau
/\tau _{0}$ and limiting cluster size distribution $\pi $, for each $\tau >0$.
\end{Prop}

We do not give the proofs of Propositions \ref{Prop:HHL1} and \ref{Prop:HHL2}
since they are readily obtained through \cite{FP} substituting Lemma 2.1  by
our Lemma \ref{Le:independencePP}. We are now prepared to give a proof of Theorem \ref{Th:theta}.
\bigskip

\begin{prooft}
\textbf{of Theorem \ref{Th:theta}} Proof of (i). First we show that (A)$\Rightarrow$(B). By Lemma \ref{Le:independencePP}, we have 
\begin{equation*}
\mathbb{P}\left( \,M_{\mathbf{W}_\rho}^{\eta }\leq v_{\rho }(\tau _{0})\,\right)
=\left( \mathbb{P}\left( \,M_{\mathfrak{Q}_{\rho}}^{\eta }\leq v_{\rho }(\tau
_{0})\,\right) \right) ^{\rho \cdot (\lambda _{d}(\mathfrak{Q}_{\rho}))^{-1}}+o(1).
\end{equation*}
Since $\lim_{\rho \rightarrow \infty }\mathbb{P}\left( \,M_{\mathbf{W}_\rho}^{\eta }\leq v_{\rho }(\tau _{0})\,\right) =e^{-\theta \tau _{0}}$ and
since $\{M_{\mathfrak{Q}_{\rho}}^{\eta }\leq v_{\rho }(\tau _{0})\}$ if and only if  $\{\#\Phi_{\mathfrak{Q}_{\rho}}^{\eta }(\tau _{0})>0\}$, it follows
that 
\begin{equation*}
\mathbb{P}\left( \,\#\Phi_{\mathfrak{Q}_{\rho}}^{\eta }(\tau _{0})>0\,\right) 
\underset{\rho \rightarrow \infty }{\sim }\frac{\lambda _{d}(\mathfrak{Q}_{\rho})}{\rho }\cdot \theta \tau _{0}.  \label{eq:limitNB0}
\end{equation*}
This together with  Proposition \ref{Prop:identification} implies that $\pi _{k}=\lim_{\rho \rightarrow \infty }\pi _{k,\mathfrak{Q}_{\rho}}(\tau_0)$
exists and  $\pi_k=p_{k}/(k\cdot \theta )$ for any $k\geq 1$. Since $\pi _{k,\mathfrak{Q}_{\rho}}(\tau _{0})\leq a_{k}$ with $\sum_{k=1}^{\infty
}a_{k}<\infty $, it follows from the dominated convergence theorem that $\pi
:=(\pi _{k})_{k\geq 1}$ is a probability measure on $\mathbf{N}_{+}$.
Applying Proposition \ref{Prop:HHL1}, we deduce that $\Phi^\eta _{\mathbf{W}_\rho}(\tau_0)$
converges to a homogeneous compound Poisson point process with intensity $\nu(\tau_0):=\theta \tau _{0}>0$
and cluster size distribution $\pi $. This together with Proposition \ref{Prop:HHL2} proves Assertion (B). 

Secondly, we show that (B)$\Rightarrow$(A). The fact that the extremal
index exists and is positive is a consequence of the fact that 
\begin{equation*}
\lim_{\rho \rightarrow \infty }\mathbb{P}\left( \,M_{\mathbf{W}_\rho}^{\eta }(1)\leq
v_{\rho }(\tau_0 )\,\right) =\lim_{\rho \rightarrow \infty }\mathbb{P}\left(
\,\#\Phi^\eta _{\mathbf{W}_\rho}(\tau_0 )=0\,\right) =e^{-\theta\tau_0},
\end{equation*}
where $\theta:=\nu(\tau_0)/\tau_0\in (0,1]$. By applying
Proposition \ref{Prop:identification}, we show that the limit of $p_k:=p_{k,\mathfrak{Q}_{\rho}}(\tau_0)$ exists and $p_k=k\theta\pi_k$. This proves Assertion (A). 

Proof of (ii). The fact that $p_k=k\theta\pi_k$ is established above. Moreover,  we have $\sum_{k=1}^\infty k^{-1}p_k = \theta\sum_{k=1}^\infty\pi_k=\theta$ since $\pi=(\pi_k)_{k\geq 1}$ is a probability measure. 
\end{prooft}

\section{Numerical illustrations}

\label{sec:examples} 

\paragraph{Layout}
In this section, we illustrate our main theorem throughout simulations for three geometric
characteristics which satisfy Condition (C). Each geometric characteristic is chosen in such a way that the value of the extremal index is known or can be conjectured. For sake of simplicity, we only do our simulations in the particular setting $d=2$. We provide approximations of  $p_1,\ldots, p_{9}$ and of the extremal index by using the fact that $\theta=\sum_{k=1}^\infty k^{-1}p_k$ (see Theorem \ref{Th:theta}, (ii)) and we compare this approximation to the theoretical value of $\theta$. 

For each geometric characteristic $g$, we proceed as follows. We take $\tau=1$ and $\rho=\exp(100)$. In particular, the cube $\mathfrak{Q}_{\rho}$, as defined in \eqref{def:B0}, is approximatively \[\mathfrak{Q}_{\rho}\simeq [-173, 173]^2,\]
by taking  $q_\rho=(\log\log\rho)^{\log\log\rho}\simeq 1134$ and $\varepsilon=0.01$. Then, we compute theoretically $v_\rho(1)$ so that $\rho\cdot \PPP{g(\cell)>v_\rho(1)}\conv[\rho]{\infty}1$. We simulate 10000 realizations of independent  Poisson-Voronoi tessellations given that the typical cell is an exceedance, i.e. $g^{\eta\cup\{0\}}(0)>v_\rho(1)$ (see Lemma \ref{Le:PalmVoronoi}). This sample of size 10000 is divided into 100 sub-samples of size 100. For each $1\leq i\leq 100$ and for each $1\leq k\leq 9$, we denote by $\hat{p}^{(i)}_k$ the empirical mean of $p_k$, i.e. the mean number of realizations in which there exist exactly $k$ Voronoi cells with nucleus in $\mathfrak{Q}_{\rho}\simeq [-173, 173]$ and such that the geometric characteristic is larger than $v_\rho(1)$. 

We summarize our empirical results by box plots associated with the empirical values $(\hat{p}^{(i)}_k)_{1\leq i\leq 100}$. For each geometric characteristic, we explain how we simulate a Poisson-Voronoi tessellation conditional on the fact that $g^{\eta\cup\{0\}}(0)>v_\rho(1)$.

\subsection{Inradius}
\label{sec:inradius}
For any $x\in \eta\subset\RR^d $, we define the so-called inradius of the Voronoi cell $C_{\eta }(x)$ as  
\begin{equation*}
r^{\eta }(x):=r(C_{\eta}(x)):=\sup \{r\geq 0:B(x,r)\subset C_{\eta }(x)\},
\end{equation*}
where $B(x,r)$ is the ball centered at $x$ with radius $r$. The
distribution of $r(\mathcal{C})$, where $r(\mathcal{C})=r(C_{\eta \cup
\{0\}}(0))$ is the typical inradius, is given by $\PPP{r(\cell)>v} = \PPP{\eta \cap B(0,2v)\neq \emptyset} = e^{-2^d\kappa_dv^d}$ for each $v\geq 0$. Hence, for any $\tau>0$, we have $\rho \cdot \mathbb{P} (r(\mathcal{C})>v_\rho(\tau) )=\tau$, when 
\begin{equation*}
v_\rho(\tau):=2^{-1}\kappa_d^{-1/d}\left(\log (\rho \tau ^{-1})\right)^{1/d}.
\end{equation*}
Moreover it is proved in \cite{CC} that \[\PPP{\max_{x\in \eta \cap \mathbf{W}_\rho}r^{\eta }(x)\leq v_\rho(\tau)}\conv[\rho]{\infty}e^{-\tau}.\]
Actually, the convergence was established for a fixed window
and for a Poisson point process such that the intensity goes to infinity. By
scaling property of the Poisson point process, the result can be re-written
as above for a fixed intensity and for a window $\mathbf{W}_\rho$ as $\rho$ goes to infinity. Therefore, we  deduce that the extremal index of the
inradius of a Poisson-Voronoi tessellation exists and is equal to $\theta=1$.
Actually, according to Theorem 2 in \cite{Chen}, the point process of exceedances $\Phi^\eta _{\mathbf{W}_\rho}(\tau)$ converges to a simple Poisson point
process of intensity $\tau$ in $W$. In particular, the distributions $\pi $ and $p$ are equal to the dirac measure at $1$.

Now, we explain how we evaluate by simulation the value of the extremal index and the
distribution $p$ when $d=2$. It is known (see e.g. \cite{Mi}) that for each $v\geq 0$, we have 
\begin{equation*}
\left( \eta \cup \{0\}|r^{\eta \cup \{0\}}(0)=v\right) \overset{\mathcal{D}}{=}\eta
_{B(0,2v)^{c}}\cup \{(2v)X_{0}\}\cup\{0\}, 
\end{equation*}
where $\eta _{B(0,2v)^{c}}$ is a Poisson point process of intensity measure $\mathbb{I}_{\{x\in B(0,2v)^{c}\}}\mathrm{d}x$ and where  $X_{0}$ is a random point uniformly distributed  on the boundary of 
$B(0,1)$. Hence, to simulate a Poisson-Voronoi tessellation provided that $r^{\eta\cup\{0\}}(0)>v_{\exp(100)}(1)\simeq 2.82$, we first simulate a random variable $r$ with distribution given by  $\PPP{r>v}=e^{-4\pi v^2}$, conditional on the fact that $r>2.82$. Then we generate a Poisson-Voronoi tessellation associated with the point process $\eta_{B(0,2r)^c}\cup \{(2r)X_{0}\}\cup\{0\}$. 

\begin{center}
\begin{figure}
\begin{center}
\begin{tabular}{cccc}
     \includegraphics[width=7.5cm,height=7.5cm]{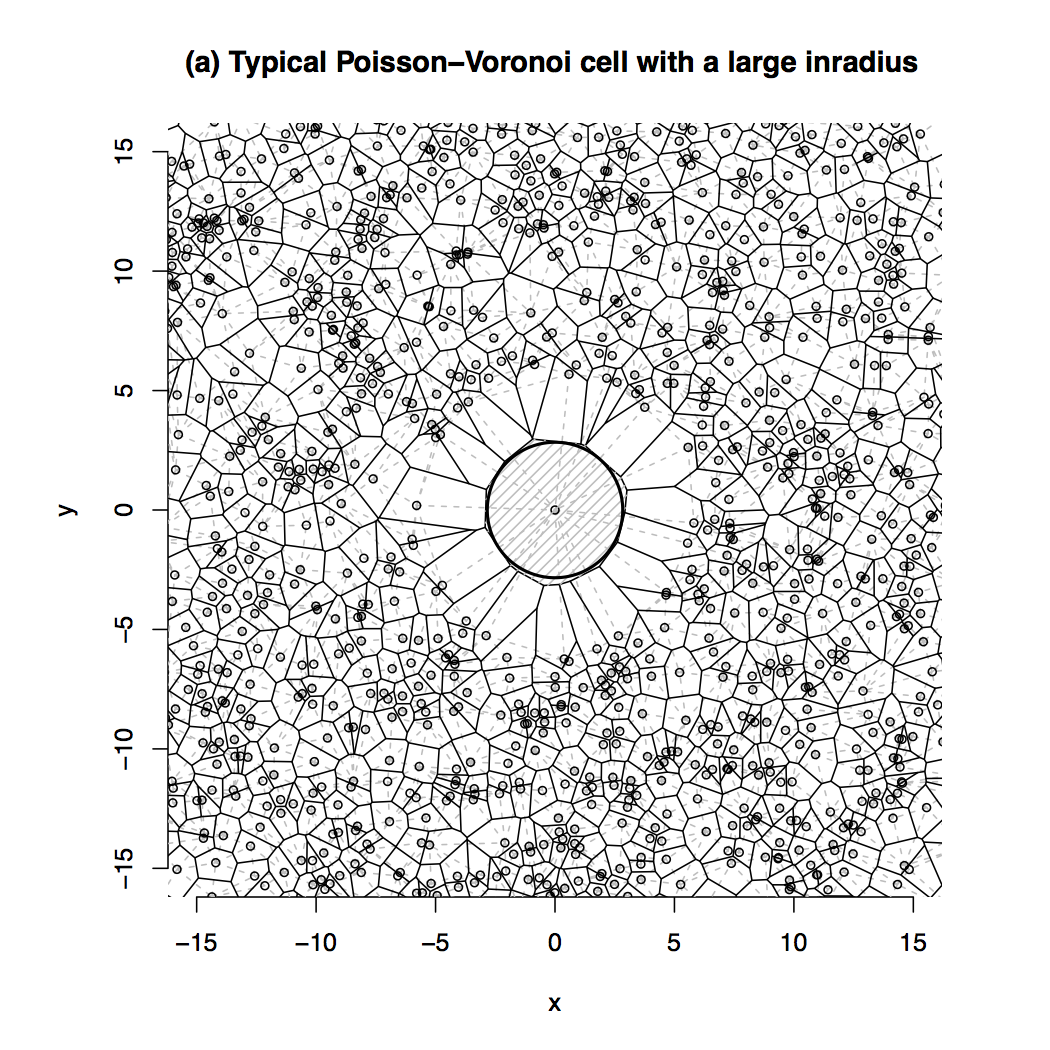} &  \includegraphics[width=7.5cm,height=7.5cm]{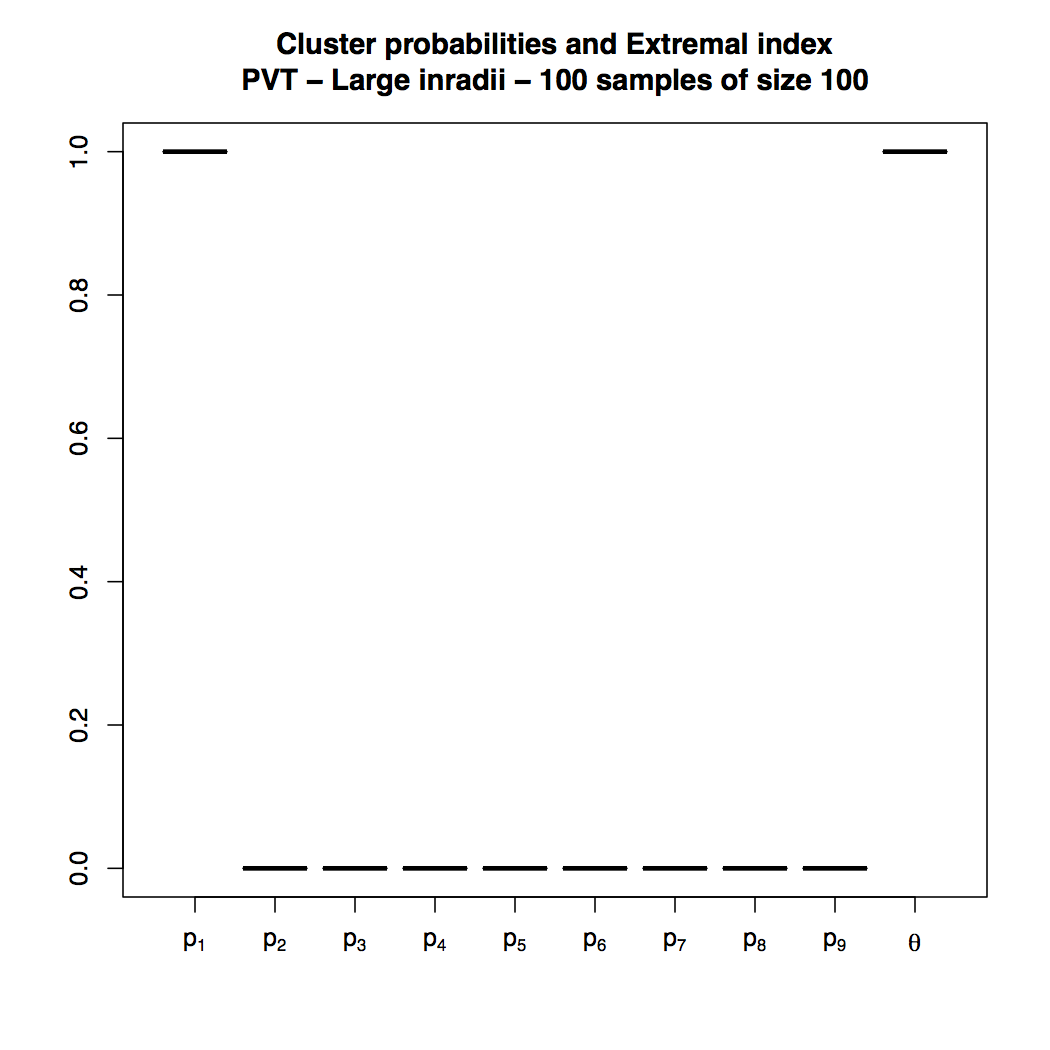}
\end{tabular}
\end{center}
\caption{\label{fig:largeinradii} Large inradius for a Poisson-Voronoi tessellation } 
\end{figure}
\end{center}

On the left part of Figure \ref{fig:largeinradii}, we provide a simulation of a Poisson-Voronoi tessellation given that $r^{\eta\cup\{0\}}(0)>2.82$. We notice that the typical cell has a shape
which tends to be circular. Actually, such an observation is related  to the D. G. Kendall's conjecture which claims that the shape of the typical Poisson-Voronoi cell in $\RR^d$, given that the volume of the cell goes to infinity, tends a.s. to a ball in $\mathbf{R}^{d}$. Many results concerning typical cells with a large geometric characteristic can be found in \cite{CalS} and \cite{HRS}. On the left part of Figure \ref{fig:largeinradii}, we also notice that there is no cell with a large inradius, excepted the typical cell. This confirms  that the cluster of exceedances are of size 1, i.e. $p_1=1$ and $\theta=1$. The right part of Figure \ref{fig:largeinradii} provides the box plots of the empirical distributions. In particular, for all simulations, we notice that there is always exactly one cell with a large inradius.  

\subsection{Reciprocal of the inradius}
\label{sec:smallinradius}

In this example, we consider the large values of the reciprocal of the inradii for a Poisson-Voronoi tessellation in $\RR^d$. Equivalently, 
this consists of the small values of the inradii. Since $\PPP{r(\mathcal{C})<v}=1-e^{-2^d\kappa_dv^d}$, we have  $\rho \cdot \mathbb{P}\left( \,r(\mathcal{C}) < v_\rho(\tau)\right)\conv[\rho]{\infty}\tau $, when 
\[v_\rho(\tau):= 2^{-1}(\kappa_d\rho)^{-1/d}\tau^{1/d} . \] Here, we have written  ``$r(\cell)< v_\rho(\tau)$'' instead of ``$r(\cell)>v_\rho(\tau)$'' in the probability because we consider the smallest inradii. Moreover, according to \cite{CC}, we know that 
\begin{equation*}
\mathbb{P}\left( \,\min_{x\in \eta \cap \mathbf{W}_\rho}r(C_\eta (x))\geq v_\rho(\tau)\right) \underset{\rho \rightarrow
\infty }{\longrightarrow }e^{-\tau /2}.
\end{equation*}
We deduce that the extremal index of the reciprocal of the inradius of a
Poisson-Voronoi tessellation exists and equals $\theta=1/2$.

\begin{center}
\begin{figure}
\begin{center}
\begin{tabular}{cccc}
    \includegraphics[width=7.5cm,height=7.5cm]{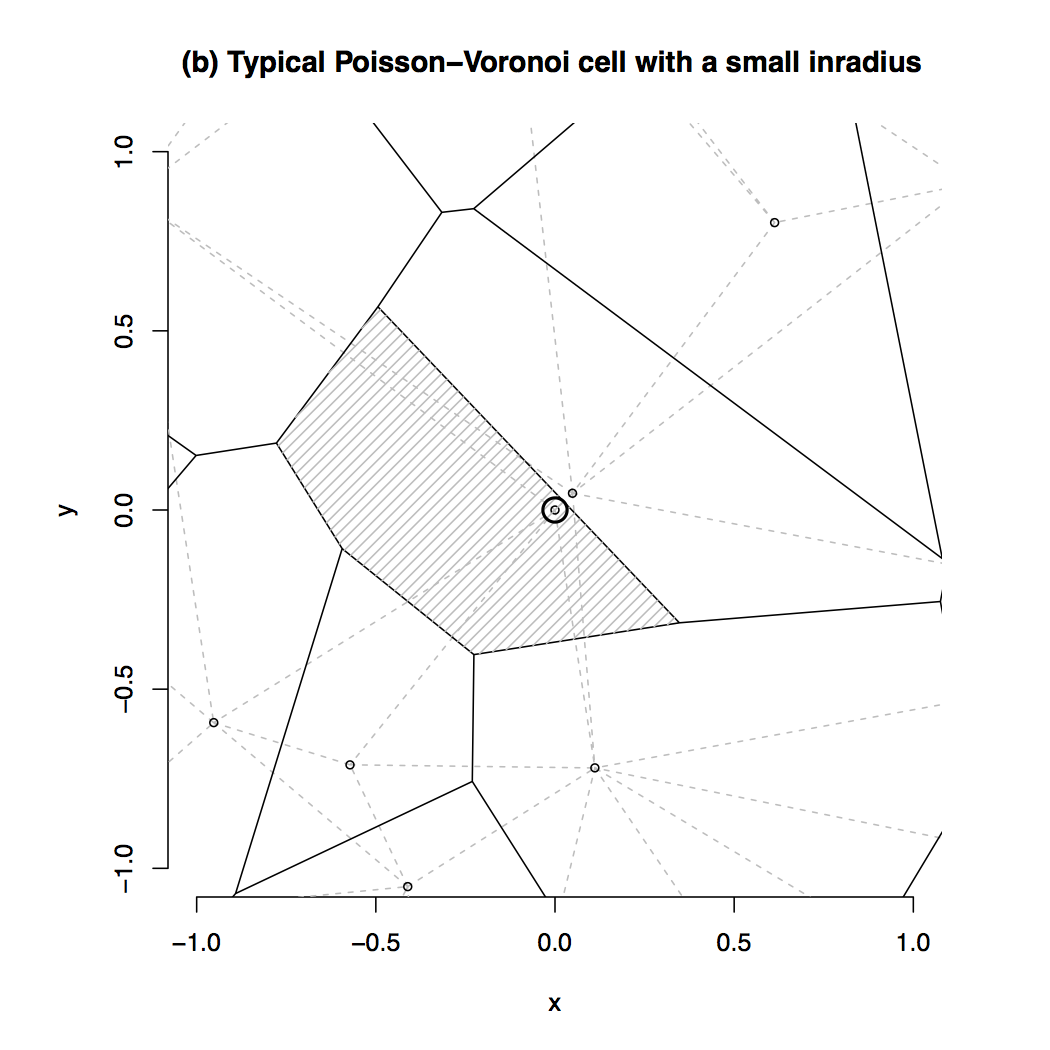} &  \includegraphics[width=7.5cm,height=7.5cm]{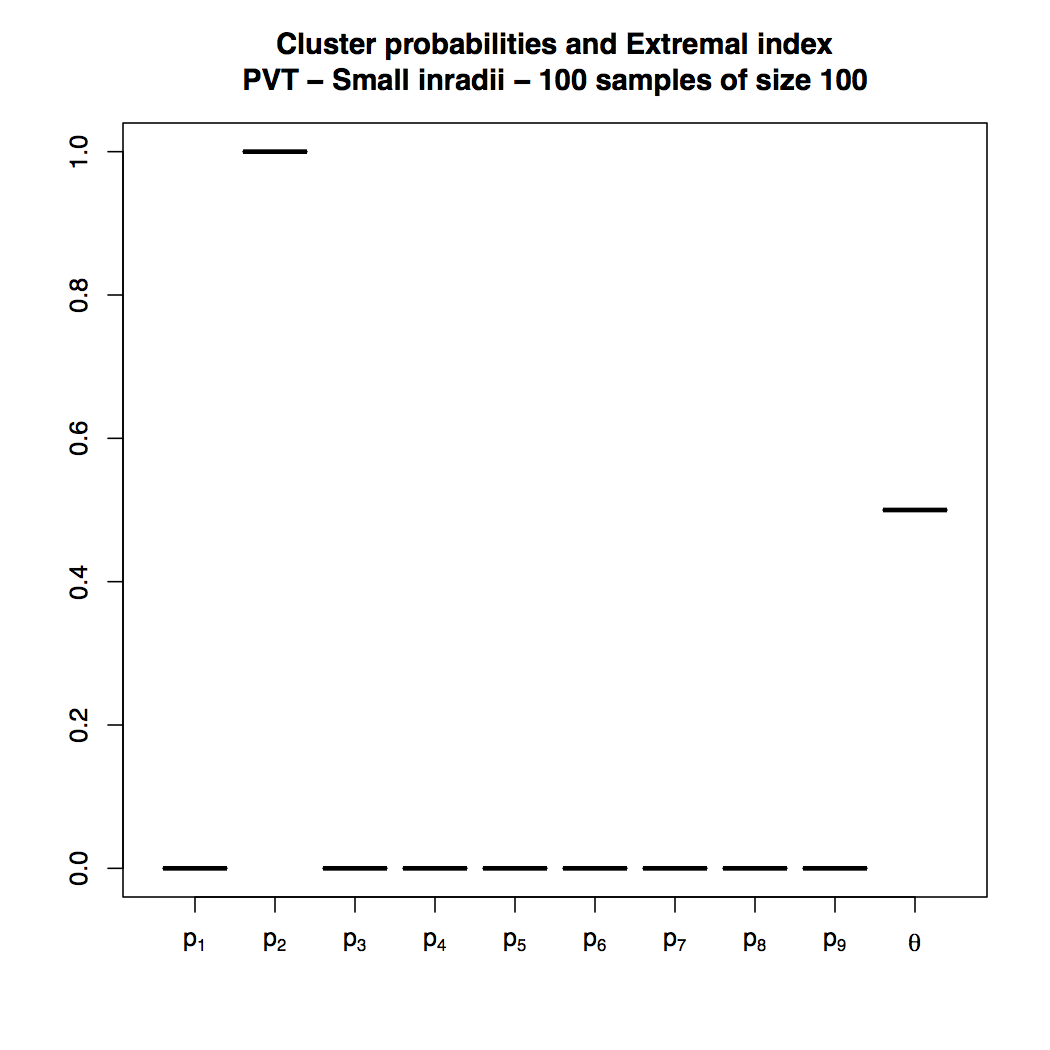}
\end{tabular}
\end{center}
\caption{\label{fig:smallinradii} Small inradius for a Poisson-Voronoi tessellation } 
\end{figure}
\end{center}

As in Section \ref{sec:inradius}, we can easily simulate a Poisson-Voronoi tessellation in $\RR^2$, conditional on the fact that $r^{\eta\cup\{0\}}(0)<v_{\exp(100)}(1)\simeq 5.44\cdot 10^{-23}$. The left part of Figure \ref{fig:smallinradii} provides a realization of a Poisson-Voronoi tessellation when $r^{\eta\cup\{0\}}(0)<v_{\exp(4)}(1)\simeq 0.0381$ (here, we have taken the threshold $v_{\exp(4)}(1)$ instead of $v_{\exp(100)}(1)$ for convenience). The fact that $\theta=1/2$ can be explained by a trivial heuristic argument: if a cell minimizes the inradius, one of its neighbors has to do the same (see also the left part of Figure \ref{fig:smallinradii}). Moreover, we can easily prove that the probability that there is more than one such a cell is negligible. Therefore clusters
are necessarily of size $2$, i.e. $p_2=1$. The right part of Figure \ref{fig:smallinradii} provides the box plots of
the empirical distributions. In particular, for all simulations, we notice that there
are always exactly two cells with a small inradius.  

\subsection{Circumradius}
For $x\in \eta\subset\RR^2 $, we define the so-called circumradius of $C_{\eta }(x)$ as
\begin{equation*}
R^\eta(x):=R(C_{\eta }(x)):=\inf \{r\geq 0:B(x,r)\supset C_{\eta }(x)\}.
\end{equation*}
According to \cite{Cal3}, we know that \[2v e^{-v }\leq \mathbb{P}\left( \,\pi R(\mathcal{C})^{2}>v
,\right) \leq 4v e^{-v },\] for each $v \geq 0.337$.  Actually, simulations suggest that the upper bound above is the order of $\mathbb{P}\left( \pi R(\mathcal{C})^{2}>v\right) $ as $v$ goes to infinity (see Table 1 in \cite{Cal3}). If we assume that $\mathbb{P}\left( \,\pi R(\mathcal{C})^{2}>v \,\right) \underset{\rho\rightarrow \infty }{\sim }a v e^{-v }$ for some $2\leq  a \leq 4$, we have $\rho \cdot \mathbb{P}\left( \, R(\mathcal{C})>v_\rho(\tau)\right)\conv[\rho]{\infty}\tau$, when 
\begin{equation*}
v_\rho(\tau):=\pi^{-1/2}\left(\log \left( a\rho \log \rho \tau ^{-1}\right)\right)^{1/2}.  
\end{equation*}
Thanks to (2.c) in \cite{CC}, we know that 
\begin{equation*}
\mathbb{P}\left( \, \max_{x\in \eta \cap \mathbf{W}_{\rho }}R^\eta
(x)\leq v_\rho(\tau) \right) \underset{\rho \rightarrow \infty }{\longrightarrow }e^{-\tau /a}.
\end{equation*}
Hence, provided that $\mathbb{P}\left( \,\pi R(\mathcal{C})^{2}>v \,\right) 
\underset{\rho \rightarrow \infty }{\sim }a v e^{-v }$, the extremal
index of the maximum of circumradius of a planar Poisson-Voronoi
tessellation exists and should be equal to $\theta=1/a$.

Now, we explain how we evaluate by simulation the value of the extremal index and the
distribution $p$. According to Lemma 1 in \cite{FZ}, we know that $R^{\eta \cup \{0\}}(0)>v$ if and only if there exists a disk of radius $v$ containing the origin on its boundary and no particle inside. Without loss of generality, we can assume that the disk, that contains the
origin on its boundary and no particle inside, has its center on the $x$-axis, since the Poisson point process is isotropic. Hence, we proceed as follows. First, we simulate a random variable $R_b$, with distribution such that $\PPP{\pi R_b^2>v}\eq[v]{\infty}bve^{-v}$, with $b=4$, given that $R_b>v_{\exp(100)}(1)\simeq 5.81$. We have taken $b=4$ since we should have $\PPP{\pi R(\cell)^2>v}\eq[v]{\infty}4ve^{-v}$ as suggested in the simulations in \cite{Cal3}. However, this choice is arbitrary and does not have influence on the final result since the conditional distribution of $R_b$ does not depend on $b$ for high thresholds. Then, we generate a Voronoi tessellation induced by the point process $\eta_{B((R,0), R)^c}\cup\{0\}$, where  $\eta_{B((R,0), R)^c}$ is a Poisson point process of intensity measure $\ind{x\in B((R,0), R)^c}\mathrm{d}x$.

\begin{center}
\begin{figure}
\begin{center}
\begin{tabular}{cccc}
   \includegraphics[width=7.5cm,height=7.5cm]{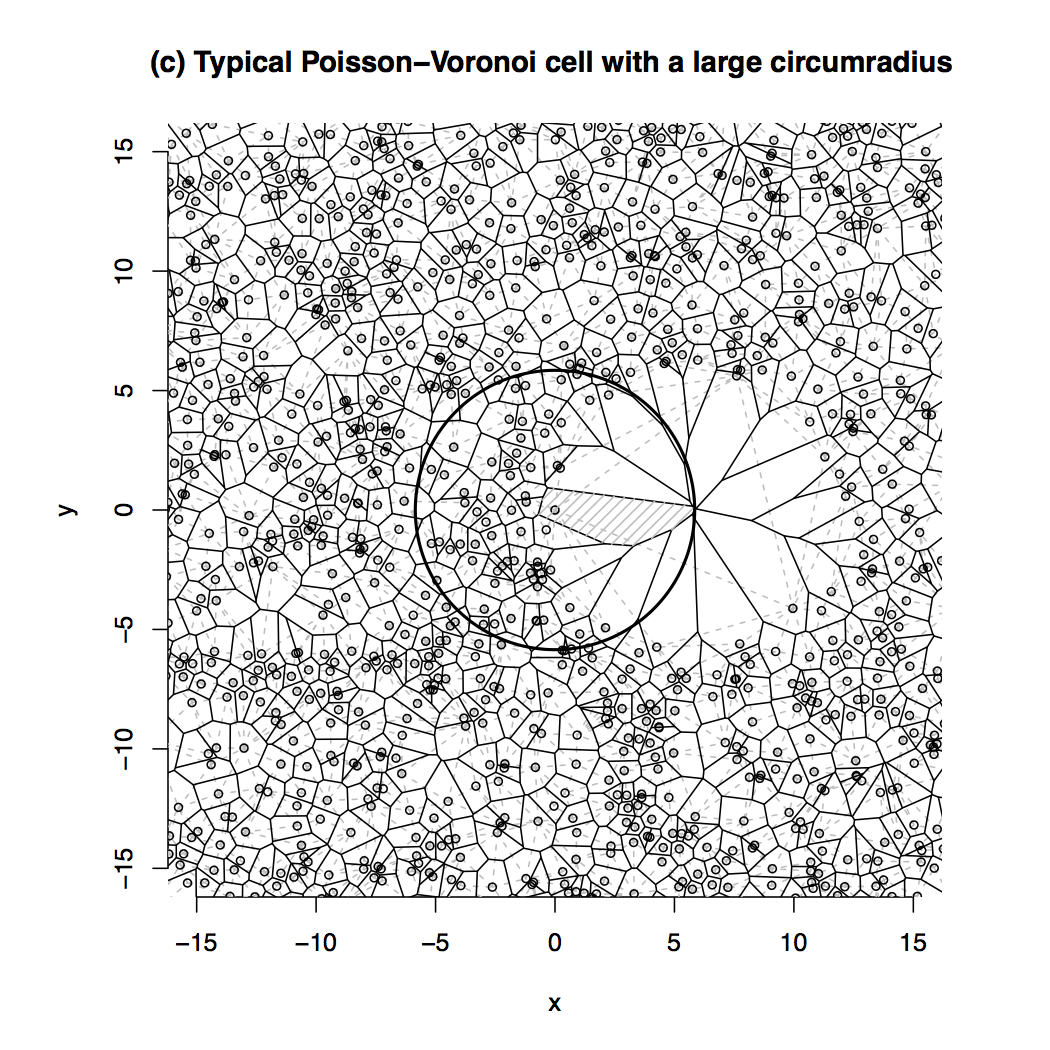} &  \includegraphics[width=7.5cm,height=7.5cm]{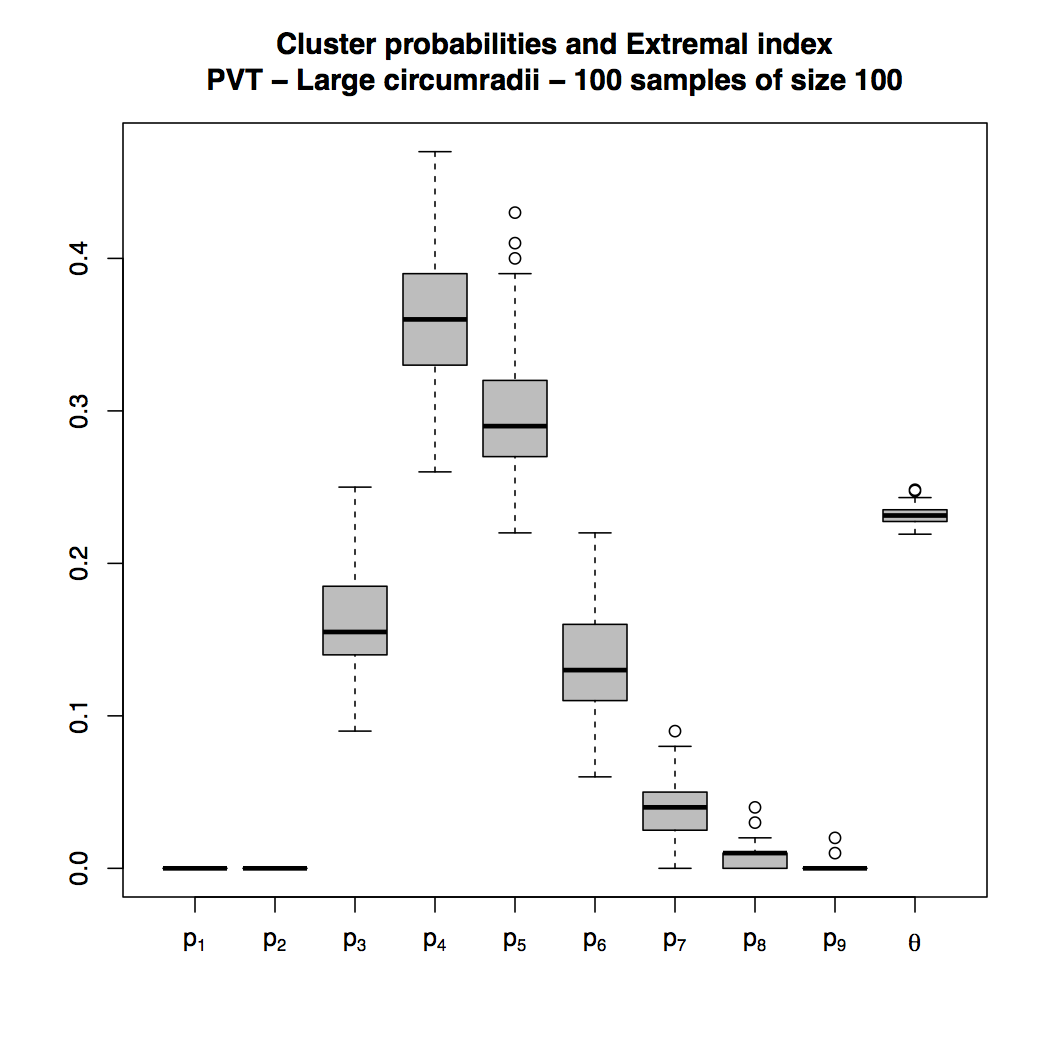}
\end{tabular}
\end{center}
\caption{\label{fig:largeCircumradii} Large circumradius for a Poisson-Voronoi tessellation } 
\end{figure}
\end{center}

On the left part of Figure \ref{fig:largeCircumradii}, we  provide a simulation of the Palm
version of the Poisson-Voronoi tessellation, given that  $R(\cell)>v_{\exp(100)}(1)\simeq 5.81$. We notice that the typical cell is very elongated and that the same fact holds for a large number of its
connected cells. In particular, the size of a cluster of exceedances is random.  On the right part of Figure \ref{fig:largeCircumradii}, we provide the box plots of the empirical distributions. This time, the
empirical distributions of the cluster size probabilities are not
degenerated for $k=3,\ldots ,9$, and their interquartile ranges are quite
large for $k=3,4,5$. We also notice that the empirical value of the
extremal index is very concentrated around a value close to $1/4$. This confirms that if $a$ exists, it should be close to 4.

\section{The case of the Poisson-Delaunay tessellation}

\label{sec:Delaunay}
\paragraph{The Poisson-Delaunay tessellation}
Let $\chi\in \mathcal{F}_{lf}$ be a locally finite subset of $\RR^d$ such that each subset of size $n<d+1$ of points are affinely independent and no $d+2$ points lie on a sphere. If two points $x,y\in\chi$ are Voronoi neighbors, i.e. $C_\chi(x)\cap C_\chi(y)\neq\varnothing$, we connect these two points by an edge. The family of these edges defines a partition of $\RR^d$ into simplices which is the so-called Delaunay tessellation. Another useful characterization of the Delaunay tessellation is the following: a simplex associated with $d+1$ points of $\chi$ is a Delaunay simplex if and only if its circumball contains no point of $\chi$ in its interior. Delaunay tessellations are very popular structures in
  computational geometry~\cite{aurenhammer2013voronoi}
and are extensively used in many areas such as surface reconstruction~\cite{cazals2006delaunay}
or mesh generation~\cite{cheng2012delaunay}. 

For each cell $C$ of the Delaunay tessellation, the nucleus $z(C)$ is defined as the center of the circumball of $C$. The set of this nuclei is denoted by $Z(\chi )$. Besides, for each $z\in Z(\chi )$, we denote by $C(z)$ the Delaunay cell whose the center of its circumball is $z$.

When $\chi =\eta $ is a homogeneous Poisson point process, the family of these 
cells is the so-called \textit{Poisson-Delaunay tessellation}. If we denote by $\gamma_\eta$ the intensity of $\eta$, then the intensity of the Poisson-Delaunay tessellation is $\gamma _{Z(\eta )}=\beta_d^{-1}\cdot \gamma_\eta$, where 
\begin{equation*}
\beta_d:=\frac{(d^{3}+d^{2})\Gamma \left( \frac{d^{2}}{2}\right) \Gamma
^{d}\left( \frac{d+1}{2}\right) }{\Gamma \left( \frac{d^{2}+1}{2}\right)
\Gamma ^{d}\left( \frac{d+2}{2}\right) 2^{d+1}\pi ^{\frac{d-1}{2}}}. 
\end{equation*}
In particular, if $d=2$, we have $\beta_2=1/2$. In the rest of the paper, we assume that $\gamma_\eta=\beta_d$ to ensure that  $\gamma _{Z(\eta )}=1$. 

The typical cell of a Poisson-Delaunay tessellation can be made explicit as follows. Let $S^{d-1}:=\{x\in \mathbf{R}^{d}:|x|=1\}$ be the unit sphere of $\mathbf{R}^{d}$ and, for $u_{1:d+1}\in (S^{d-1})^{d+1}$, let $\Delta(u_{1:d+1}):=conv(u_{1},\ldots ,u_{d+1})$. According to Miles (see e.g. Theorem 10.4.4 in \cite{SW}), for any bounded measurable function $f:\mathbf{R}_{+}\times S^{d-1}\rightarrow 
\mathbf{R}$, we have 
\begin{equation}
\mathbb{E}[f(\mathcal{C})]:=a_{d}\gamma _{\eta }^{d}\int_{\mathbf{R}_{+}}\int_{S^{d-1}}\cdots \int_{S^{d-1}}a(u_{1:d+1})r^{d^{2}-1}e^{-\gamma _{\eta
}\kappa _{d}r^{d}}f(\Delta (ru_{1:d+1}))\sigma (\mathrm{d}u_{1:d+1})\mathrm{d}r,
\label{def:TypicalDelaunay}
\end{equation}
where $a_{d}:=\beta _{d}/(d+1)$ and $a(u_{1:d+1}):=\lambda _{d}\left( \Delta
(u_{1:d+1})\right)$. The measure $\sigma(\mathrm{d}u)$ is the uniform distribution on $S^{d-1}$ with normalization $\sigma(S^{d-1})=\omega_{d-1}$, where $\omega_{d-1}:=d\kappa_d$ is the area of the unit sphere and $\sigma (\mathrm{d}u_{1:d+1}):=\bigotimes_{i=1}^{d+1}\sigma(\mathrm{d}u_i)$. Hence, the typical satisifes the equality in distribution  $\cell \overset{\mathcal{D}}{=} \Delta(RU_{1:d+1})$, where $R>0$ and $U_{1:d+1}\in(S^{d-1})^{d+1}$ are two independent random variables whose the distributions are provided in \eqref{def:typicalcell}.

\paragraph{The extremes of the Poisson-Delaunay tessellation}

Let $g$ be a geometric characteristic such that \eqref{Threshcond} holds.  As for a Poisson-Voronoi tessellation, we consider  the point process of normalized exceedances, say 
\begin{equation*}
\Psi^{\eta }(\tau):=\rho ^{-1/d}\cdot  \left\{z\in Z(\eta ):g(C(z))>v_{\rho
}(\tau )\right\}.
\end{equation*}
For any Borel subset $B\subset\RR^d$, we  write $\Psi^{\eta }_B(\tau):=\Psi^\eta(\tau)\cap (\rho^{-1/d}B)$. We also let $\Psi ^{\eta, 0 }(\tau)$ be the Palm version of $\Psi ^{\eta }(\tau)$ and  $\Psi _{B}^{\eta, 0 }(\tau)=\Psi ^{\eta, 0 }(\tau)\cap B$. In the rest of the paper, the quantity $p_{k,B}(\tau)$ refers to as the probability that there exist $k$ exceedance cells in $B$ conditional on the fact that the typical cell is an exceedance, i.e. $p_{k,B}(\tau):=\PPP{\# \psi_{B}^{\eta,0}(\tau)=k}$. In the same spirit as Lemma  \ref{Le:PalmVoronoi}, we provide below an explicit characterization of this  probability. 

\begin{Prop}
\label{Prop:DelaunayPalm}
Let $\mathcal{A}$ be a Borel subset in $\mathcal{F}_{lf}$. Then 
\begin{equation*}
\mathbb{P}\left( \,\Psi ^{\eta, 0 }(\tau)\in \mathcal{A}\,\right) =\mathbb{P}\left( \,\left.
\Psi ^{\eta _{\mathbf{R}^{d}\setminus B(0,R)}\cup \{RU_{1:d+1}\}}(\tau)\in
\mathcal{A}\right\vert g(\Delta (RU_{1:d+1}))>v_{\rho }(\tau )\right) .
\end{equation*}
Therefore, for any $B\subset \mathbf{R}^{d}$, 
\begin{equation*}
p_{k,B}(\tau)=\mathbb{P}\left( \,\left. \#\Psi _{B}^{\eta _{\mathbf{R}^{d}\setminus B(0,R)}\cup \{RU_{1:d+1}\}}(\tau)=k\right\vert g(\Delta
(RU_{1:d+1}))>v_{\rho }(\tau )\right) .
\end{equation*}
\end{Prop}


\begin{prooft}
Let $A\subset\RR^d$ be such that $\lambda_d(A)=1$. It follows from the definition of the Palm distribution of the point process $\Psi ^{\eta}(\tau)$, that 
\begin{equation*}
\begin{split}
\mathbb{P}\left( \Psi ^{\eta, 0 }(\tau)\in \mathcal{A}\,\right) & :=\frac{1}{\gamma _{\Psi
^{\eta }(\tau)}}\mathbb{E}\left[ \,\sum_{z\in \Psi ^{\eta }(\tau)\cap A}\mathbb{I}_{(\Psi ^{\eta }(\tau)-z)\in \mathcal{A}}\right]  \\
& =\frac{1}{\gamma _{\Psi^{\eta }(\tau)}}\mathbb{E}\left[
\,\sum_{\{x_{1:d+1}\}\subset \eta }\mathbb{I}_{(\Psi^{\eta }(\tau)-\rho
^{1/d}z(x_{1:d+1}))\in \mathcal{A}}\,\mathbb{I}_{z(x_{1:d+1})\in Z(\eta )\cap \mathbf{A}_\rho
}\,\,\right],
\end{split}
\end{equation*}
where $z(x_{1:d+1})$ is the center of the circumball of the simplex $\Delta(x_{1:d+1})$. According to the Slivnyak-Mecke formula and the Blaschke-Petkantschin formula (e.g. Theorem 7.3.1 in \cite{SW}) that  
\begin{multline*}
\mathbb{P}\left( \Psi ^{\eta, 0 }(\tau)\in \mathcal{A}\right) =\frac{\gamma _{\eta }^{d+1}d!}{\gamma _{\Psi ^{\eta }(\tau)}(d+1)!}\int_{\mathbf{A}_\rho}\int_{\mathbf{R}_{+}}\int_{S^{d+1}}r^{d^{2}-1}a(u_{1:d+1}) \\
\times \mathbb{P}\left( (\Psi ^{\eta }(\tau)-\rho ^{1/d}z)\in \mathcal{A}, \eta \cap
B(z,r)=\varnothing \,\right) \mathbb{I}_{g(z+\Delta (ru_{1:d+1}))>v_{\rho}(\tau)}\,\sigma(\mathrm{d}u_{1:d+1})\mathrm{d}r\mathrm{d}z.
\end{multline*}
Since $\eta $ is stationary and since $g$ is translation-invariant, the
integrand does not depend on $z$. Integrating over $z\in \mathbf{A}_\rho$, and using the fact
that $\lambda _{d}(\mathbf{A}_\rho)=\rho$, we get 
\begin{multline*}
\mathbb{P}\left( \,\Psi ^{\eta, 0 }(\tau)\in \mathcal{A}\right) =\frac{\gamma _{\eta
}^{d+1}\rho }{\gamma _{\Psi ^{\eta }(\tau)}(d+1)}\int_{\mathbf{R}_{+}}\int_{S^{d+1}}r^{d^{2}-1}a(u_{1:d+1}) \\
\times \mathbb{P}\left( \Psi ^{\eta }(\tau)\in \mathcal{A}, \eta \cap B(0,r)=\varnothing
\,\right) \mathbb{I}_{g(\Delta (ru_{1:d+1}))>v_{\rho}(\tau)}\,\sigma(\mathrm{d}u_{1:d+1})\mathrm{d}r.
\end{multline*}
We give below an explicit representation for the integrand. Let $\eta _{B(0,r)}$ and $\eta _{\mathbf{R}^{d}\setminus B(0,r)}$ be two independent Poisson point processes
with intensity measures $\gamma_\eta\ind{x\in B(0,R)}\mathrm{d}x$ and $\gamma_\eta\ind{x\in \mathbf{R}^{d}\setminus B(0,r)}\mathrm{d}x$ respectively. We know that 
\begin{equation*}
\eta \overset{\mathcal{D}}{=}\eta _{B(0,r)}\cup \eta _{\mathbf{R}^{d}\setminus B(0,r)}.
\end{equation*}
This gives 
\begin{equation*}
\begin{split}
\mathbb{P}\left( \Psi ^{\eta }(\tau)\in \mathcal{A}, \eta \cap B(0,r)=\varnothing \,\right) &
=\mathbb{P}\left( \,\Psi ^{\eta _{\mathbf{R}^{d}\setminus B(0,r)}\cup
\{ru_{1:d+1}\}}(\tau)\in \mathcal{A}, \eta _{B(0,r)}\cap B(0,r)=\varnothing \,\right)  \\
& =e^{-\gamma _{\eta }\kappa _{d}r^{d}}\mathbb{P}\left( \,\Psi ^{\eta _{\mathbf{R}^{d}\setminus B(0,r)}\cup \{ru_{1:d+1}\}}(\tau)\in \mathcal{A}\,\right) .
\end{split}
\end{equation*}
Hence,  
\begin{multline*}
\mathbb{P}\left( \,\Psi ^{\eta, 0 }(\tau)\in \mathcal{A}\right) =\frac{\gamma _{\eta
}^{d+1}\rho }{\gamma _{\Psi ^{\eta }(\tau)}(d+1)}\int_{\mathbf{R}_{+}}\int_{S^{d+1}}r^{d^{2}-1}a(u_{1:d+1}) \\
\times e^{-\gamma _{\eta }\kappa _{d}r^{d}}\mathbb{P}\left( \Psi ^{\eta _{\mathbf{R}^{d}\setminus B(0,r)}\cup \{ru_{1:d+1}\}}(\tau)\in \mathcal{A}\,\right)
\sigma(\mathrm{d}u_{1:d+1})\mathrm{d}r.
\end{multline*}
Moreover, we know that $\gamma _{\Psi ^{\eta }(\tau)}=\rho \mathbb{P}\left( g(\mathcal{C})>v_{\rho }(\tau )\right) $ and that $\gamma _{\eta }=(d+1)a_{d}$. Then, we get 
\begin{multline*}
\mathbb{P}\left( \,\Psi ^{\eta, 0 }(\tau)\in \mathcal{A}\right) =\frac{a_{d}\gamma _{\eta
}^{d}}{\mathbb{P}\left( \,g(\mathcal{C})>v_{\rho }(\tau )\right) }\int_{\mathbf{R}_{+}}\int_{S^{d+1}}r^{d^{2}-1}a(u_{1:\mathrm{d}+1}) \\
\times e^{-\gamma _{\eta }\kappa _{d}r^{d}}\mathbb{P}\left( \,\Psi ^{\eta _{\mathbf{R}^{d}\setminus B(0,r)}\cup \{ru_{1:d+1}\}}\in \mathcal{A}\,\right)
\sigma(\mathrm{d}u_{1:d+1})\mathrm{d}r.
\end{multline*}
This proves the first equality in Proposition \ref{Prop:DelaunayPalm} since $\mathcal{C}\overset{\mathcal{D}}{=}\Delta (RU_{1:d+1})$. The second equality is a direct consequence of the first one. 
\end{prooft}

We think that Theorem \ref{Th:theta} can be adapted in the context of a Poisson-Delaunay tessellation. To do it, we have to replace the point process  $\Phi ^{\eta }(\tau)$ by the point process $\Psi ^{\eta }(\tau)$ and we have to use the characterization of the probability $p_{k, \mathfrak{Q}_{\rho}}(\tau_0)$ as described in the above proposition. We can easily extend Lemma \ref{Le:Arho} and adapt Condition (C) in the particular setting of a Poisson-Delaunay tessellation. However, the main difficulty to adapt Theorem \ref{Th:theta} focuses on an analogous version of Proposition \ref{Prop:identification} since its proof  seems very technical. We give below a numerical illustration which confirms that Theorem \ref{Th:theta} should be true for a Poisson-Delaunay tessellation.

\paragraph{A numerical illustration}

Let $m_{PDT}$ be a Poisson-Delaunay tessellation generated by a 
Poisson point process $\eta$ in $\RR^d$ with intensity $\gamma_\eta=\beta_d$. For each cell $C\in m_{PDT}$, we
consider the so-called circumcenter of $C$ defined as 
\begin{equation*}
R(C):=\inf \{R\geq 0:C\subset B(z,R),z\in \mathbb{R}^{d}\}.
\end{equation*}
According to $\eqref{def:TypicalDelaunay}$, the random variable $\kappa_dR(\mathcal{C})^{d}$ is Gamma 
distributed with parameters $(d,\beta_d)$. A Taylor
expansion of $\mathbb{P}\left( R(\mathcal{C})>v \right) $ as $v$ goes
to infinity (e.g. Equation (3.14) in \cite{Chen}), shows that $\rho \cdot \mathbb{P}\left( \,R(\mathcal{C})>v_\rho(\tau)\right)\conv[\rho]{\infty}\tau$, when 
\begin{equation*}
v_\rho(\tau):= (\kappa_d\beta
_{d})^{-1/d}\cdot \left(\log \left( [(d-1)!]^{-1}\rho \log (\beta _{d}\rho
)^{d-1}\tau ^{-1}\right)\right)^{1/d}.
\end{equation*}
Moreover, with standard arguments, we can easily show that the maximum of
circumradii of Delaunay cells $\max_{\underset{}{C\in m_{PDT},z(C)\in 
\mathbf{W}_\rho}}R(C)$ has the same asymptotic behavior as the maximum of
circumradii of the associated Voronoi cells $\max_{x\in \eta \cap \mathbf{W}_\rho}R(C_{\eta }(x))$. Besides, according to (2c) in \cite{Chen}, we know that 
\begin{equation*}\PPP{\max_{x\in \eta\cap \mathbf{W}_{\rho}}R(C_{\eta}(x)) \leq (\kappa_d\beta_d)^{-1/d} \left(\log\left(\alpha_d\beta_d\rho\log(\beta_d\rho)^{d-1}\tau^{-1}\right) \right)^{1/d}}\conv[\rho]{\infty}e^{-\tau},\end{equation*} where $\alpha_d:=\frac{1}{d!}\left(\frac{\pi^{1/2}\Gamma\left(\frac{d}{2}+1\right)}{\Gamma\left(\frac{d+1}{2}\right)}\right)^{d-1}$. It follows that 
\begin{equation*}
\mathbb{P}\left( \,\max_{\underset{z(C)\in \mathbf{W}_{\rho }}{C\in m_{PDT},}}R(C)\leq v_\rho(\tau) \,\right) 
\underset{\rho \rightarrow \infty }{\longrightarrow }e^{-\theta _{d}\tau },
\end{equation*}
where 
\begin{equation*}
\theta _{d}:=\alpha _{d}\beta _{d}(d-1)!=\frac{(d^{3}+d^{2})\Gamma \left( 
\frac{d^{2}}{2}\right) \Gamma \left( \frac{d+1}{2}\right) }{d\Gamma
\left( \frac{d^{2}+1}{2}\right) \Gamma \left( \frac{d+2}{2}\right) 2^{d+1}}.
\end{equation*}
In particular, when $d=1,2,3$, the extremal index equals $\theta_1 =1$, $\theta_2 =1/2$ and $\theta_3 =35/128$ respectively.

Now, we explain how we evaluate by simulation the value of the extremal index and the
distribution $p$ when $d=2$. First, we simulate a random variable $R$ such that $\pi R^2$ is Gamma
distributed with parameters $(2,1/2)$, given that $R>v_{\exp(100)}(1)\simeq 8.16$. Then we simulate a  typical cell $\cell$, with circumradius $R$, by using the method described in \cite{Ke}. The Poisson-Delaunay tessellation which is generated is induced by the point process $\eta_{B(0,2R)^c}\cup\{0\}$, where $\eta_{B(0,2R)^c}$ is a Poisson point process with intensity measure $\ind{x\in B(0,2R)^c}\mathrm{d}x$ (see Proposition \ref{Prop:DelaunayPalm}). 

\begin{center}
\begin{figure}
\begin{center}
\begin{tabular}{cccc}
   \includegraphics[width=7.5cm,height=7.5cm]{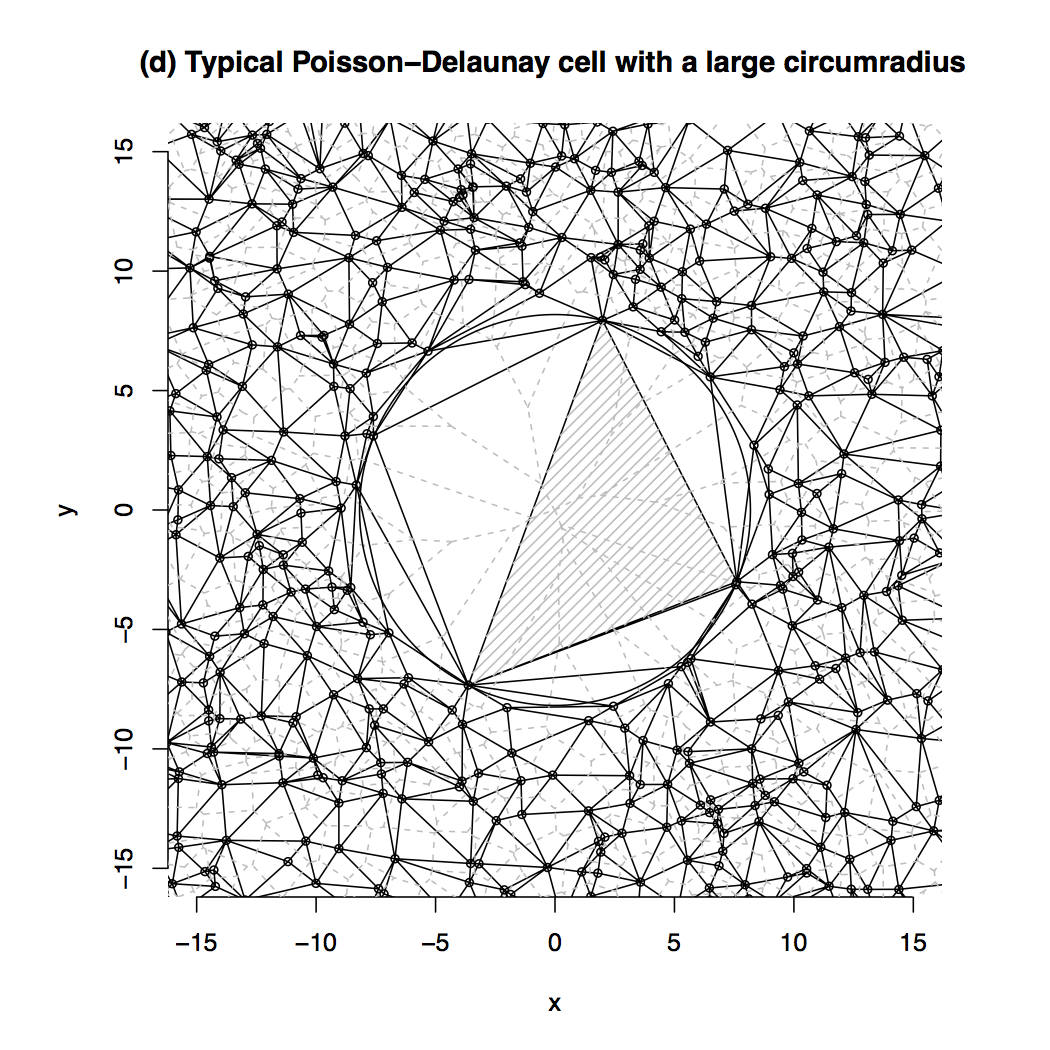} &  \includegraphics[width=7.5cm,height=7.5cm]{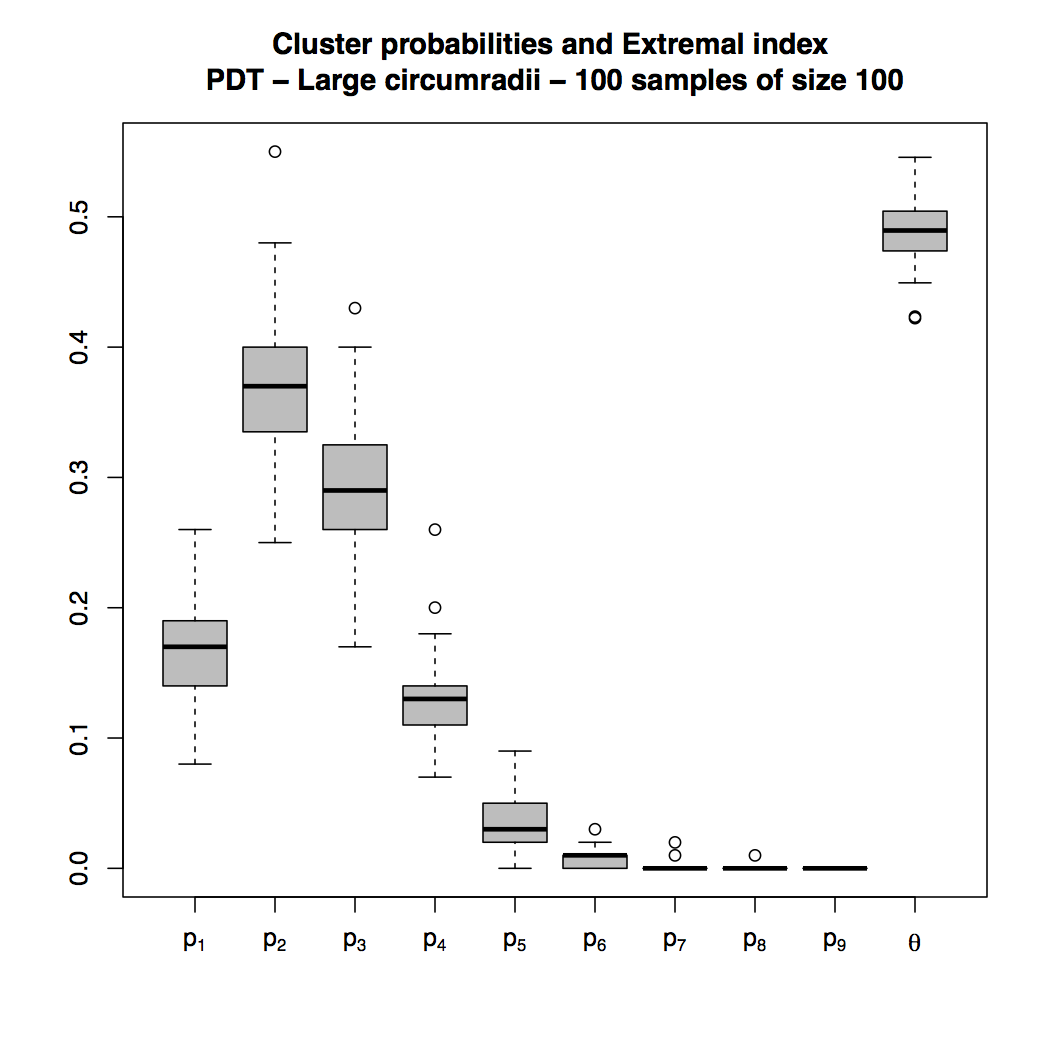}
\end{tabular}
\end{center}
\caption{\label{fig:largecircumradiidelaunay} Large circumradius for a Poisson-Delaunay tessellation } 
\end{figure}
\end{center}

On the left part of Figure \ref{fig:largecircumradiidelaunay}, we provide a simulation of the Palm version of the Poisson-Delaunay tessellation given that the typical cell has a circumradius
larger than $8.16$. The number of neighbors of the typical cells which are exceedances is random. The right part of Figure \ref{fig:largecircumradiidelaunay} provides the box plots of
the empirical probabilities. Notice that these empirical distributions  are not degenerated for $k=1,\ldots ,8$. Their interquartile ranges are not so important as for the circumradii of the Poisson-Voronoi tessellation, but the spread of the empirical distribution of the extremal index is larger. Besides, the empirical value of the extremal index is very concentrated around a value close to $1/2$, which is the theoretical value of $\theta$. 

\end{document}